\RequirePackage{fix-cm}
\documentclass{svjour3}

\smartqed 

\usepackage{graphicx}
\usepackage{float}
\usepackage{latexsym}

\usepackage{subcaption}

\usepackage[dvipsnames]{xcolor}
\usepackage{enumitem,graphicx}

\usepackage[pagebackref=true]{hyperref}
\hypersetup{
    pdftoolbar=true,        
    pdfmenubar=true,        
    pdffitwindow=false,     
    pdfstartview={FitH},    
    colorlinks=true,       
    linkcolor=OliveGreen,          
    citecolor=blue,        
    filecolor=black,      
    urlcolor=red           
}

\def\makeheadbox{\relax}

\begin{document}
\title{Can the Elliptic Billiard Still Surprise Us?}
\author{Dan Reznik \and Ronaldo Garcia \and Jair Koiller}

\authorrunning{Reznik, Garcia and Koiller}
\institute{D. Reznik \at Data Science Consulting \email{dreznik@gmail.com} \and
R. Garcia \at Universidade Federal de Goiás \email{ragarcia@ufg.br} \and 
J. Koiller \at Universidade Federal de Juiz de Fora \email{jairkoiller@gmail.com}}

\date{November 2020}

\maketitle

\begin{abstract}
Can any secrets still be shed by that much studied, uniquely integrable, Elliptic Billiard? Starting by examining the family of 3-periodic trajectories and the loci of their Triangular Centers, one obtains a beautiful and variegated gallery of curves: ellipses, quartics, sextics, circles, and even a stationary point. Secondly, one notices this family conserves an intriguing ratio: Inradius-to-Circumradius. In turn this implies three conservation corollaries: (i) the sum of bounce angle cosines, (ii) the product of excentral cosines, and (iii) the ratio of excentral-to-orbit areas. Monge's Orthoptic Circle's close relation to 4-periodic Billiard trajectories is well-known. Its geometry provided clues with which to generalize 3-periodic invariants to trajectories of an arbitrary number of edges. This was quite unexpected. Indeed, the Elliptic Billiard did surprise us!

\keywords{elliptic billiard, periodic trajectories, integrability, triangle center, locus, loci, conservation, invariance, invariant, constant of motion}

\noindent \textbf{MSC} {37-40, 51N20, 51M04, 51-04}
\end{abstract}

\section{Introduction}
\label{sec:intro}
\begin{quote}
Things in motion sooner catch the eye\\
Than what not stirs.\\
--W.S., ``Troilus and Cressida'' 
\end{quote}

\begin{figure}
    \centering
    \includegraphics[width=\textwidth]{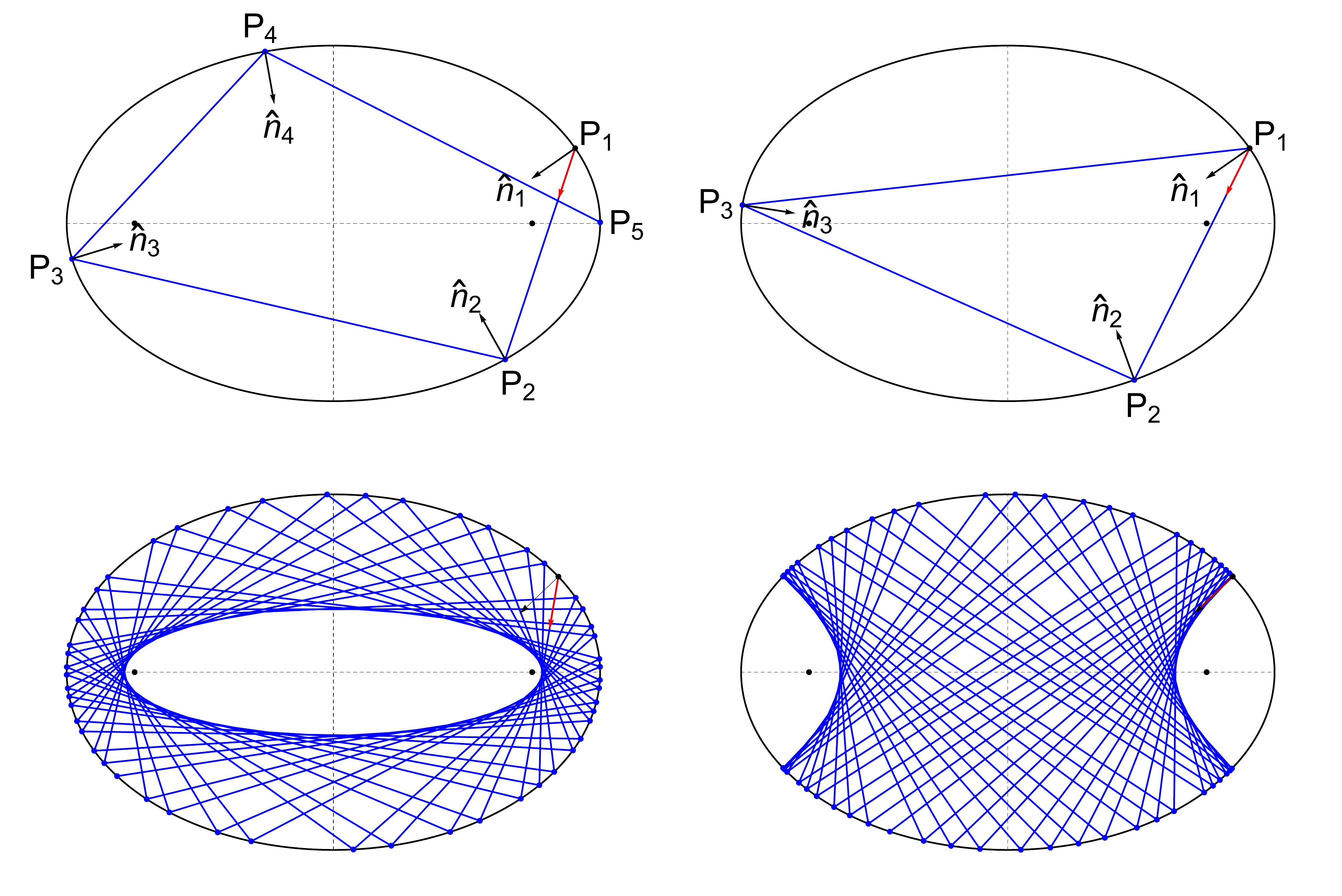}
    \caption{Particle trajectories in an Elliptic Billiard. \textbf{Top left}: the first four segments of a trajectory starting at $P_1$ along the red arrow, bouncing elastically at $P_2$ (reflecting about normals $\hat{n}_i$), etc. \textbf{Top right}: a 3-periodic trajectory. \textbf{Bottom}: quasi-periodic regimes tangent to an elliptic (left) or hyperbolic (right) {\em Caustic}. The first (resp. second) regime is established if a segment passes outside (resp. between) the foci (black dots).
    \href{https://youtu.be/A7mPzrNJHkA}{Video} \cite[pl\#1]{dsr_math_intell_playlist}}
    \label{fig:billiard-trajectories}
\end{figure}

The Elliptic Billiard (EB) consists of a particle moving with constant velocity in the interior of an ellipse, bouncing elastically against its boundary. A few trajectory regimes are shown in Figure~\ref{fig:billiard-trajectories}. Though the EB has been exhaustively studied \cite{dragovic11,sergei91}, can it still yield any surprises?

\begin{figure}
    \centering
    \includegraphics[width=.7\textwidth]{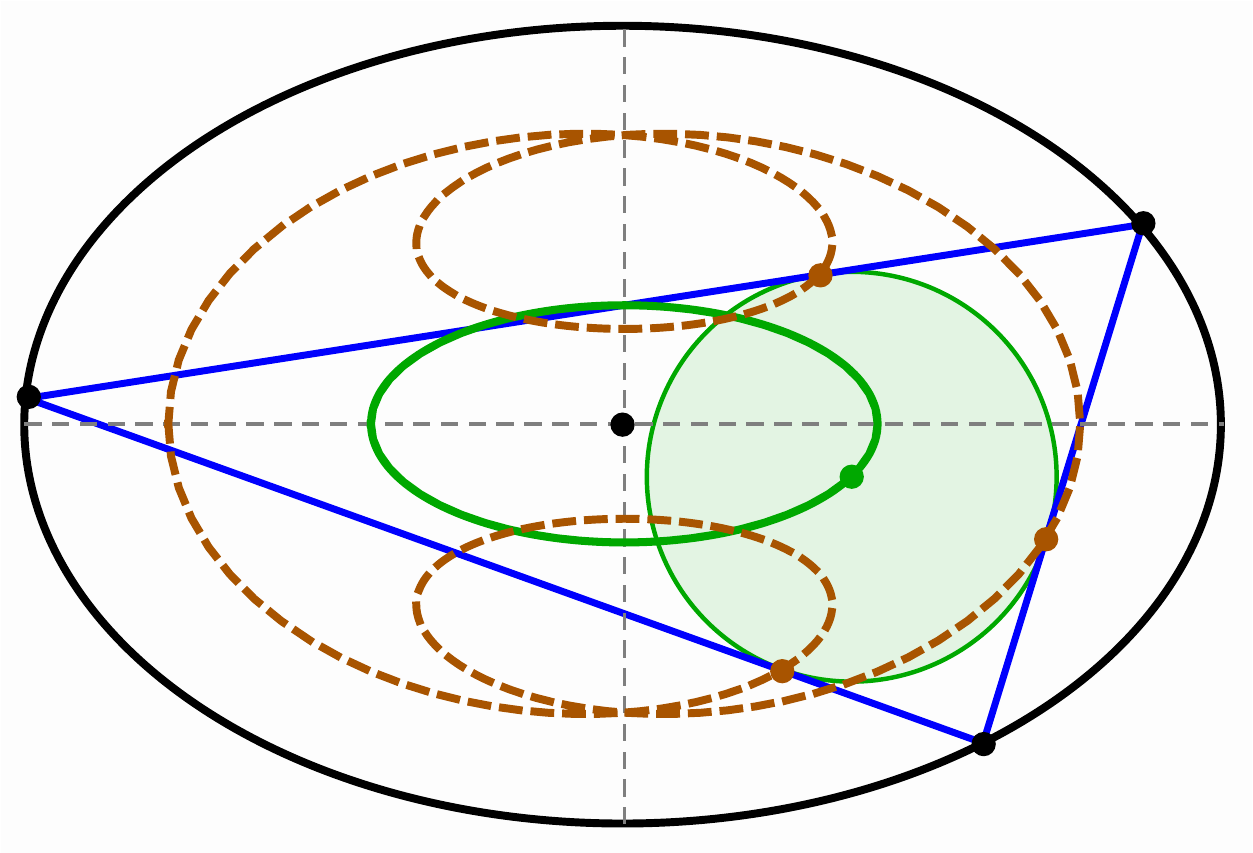}
    \caption{An $N=3$ orbit (blue), its Incircle (transparent green), Incenter (green dot) and Intouch Points (brown dots). Over the $N=3$ family, the Incenter locus is a perfect ellipse (green), while the Intouchpoints produce a self-intersecting sextic (dashed brown).
    \href{https://youtu.be/9xU6T7hQMzs}{Video} \cite[pl\#2]{dsr_math_intell_playlist}}
    \label{fig:intro-plot}
\end{figure}

\begin{figure}
    \centering
    \includegraphics[width=\textwidth]{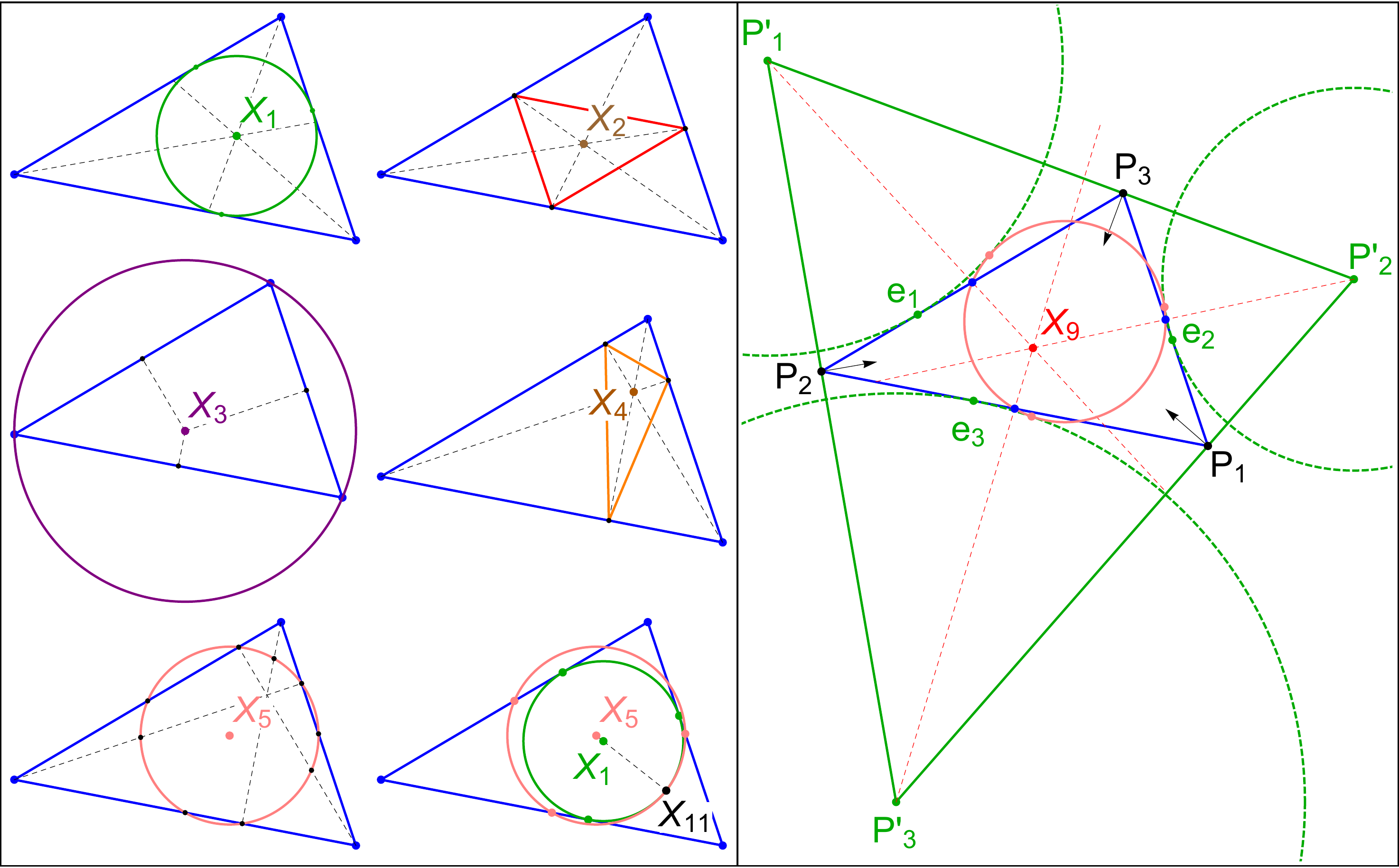}
    \caption{Notable Triangle Points, referred to as $X_i$, after Kimberling \cite{etc}. \textbf{Left}: The Incenter $X_1$ is the intersection of angular bisectors, and center of the Incircle (green), a circle tangent to the sides at three {\em Intouchpoints} (green dots), its radius is the {\em Inradius} $r$. The Barycenter $X_2$ is where lines drawn from the vertices to opposite sides' midpoints meet. Side midpoints define the {\em Medial Triangle} (red). The Circumcenter $X_3$ is the intersection of perpendicular bisectors, the center of the {\em Circumcircle} (purple) whose radius is the {\em Circumradius} $R$. The Orthocenter $X_4$ is where altitudes concur. Their feet define the {\em Orthic Triangle} (orange). $X_5$ is the center of the 9-Point (or Euler) Circle (pink): it passes through each side's midpoint, altitude feet, and Euler Points \cite{mw}. The Feuerbach Point $X_{11}$ is the single point of contact between the Incircle and the 9-Point Circle. \textbf{Right}: given a reference triangle $P_1P_2P_3$ (blue), the {\em Excenters} $P_1'P_2'P_3'$ are pairwise intersections of lines through the $P_i$ and perpendicular to the bisectors. This triad defines the {\em Excentral Triangle} (green). The {\em Excircles}  (dashed green) are centered on the Excenters and are tangent to the triangle sides at the {\em Extouch Points} $e_i,i=1,2,3$. Lines drawn from each Excenter through sides' midpoints (dashed red) concur at the {\em Mittenpunkt} $X_9$. Also shown is Feuerbach's Theorem \cite{mw}: besides the Incircle, the 9-Pt Circle (pink) touches each Excircle at a single point (pink dots), the vertices of the {\em Feuerbach Triangle}.}
    \label{fig:constructions}
\end{figure}

In 2011, we uploaded a
\href{https://www.youtube.com/watch?v=9xU6T7hQMzs}{video} \cite[pl\#2]{dsr_math_intell_playlist} showing the family of 3-periodic trajectories as a continuous set of rotating triangles. We also drew the locus of their Incenter and {\em Intouchpoints}\footnote{Where the Incircle touches a triangle's sides}. Though the latter was a self-intersecting sextic, Figure~\ref{fig:intro-plot}, we were surprised the former was (numerically) a perfect ellipse. For a refresher on triangle concepts see Figure~\ref{fig:constructions}. 

 Over the next few years there surfaced elegant proofs showing that the locus of the Incenter \cite{ronaldo16,olga14}, Barycenter \cite{ronaldo19,sergei2016}, Circumcenter  \cite{corentin19,ronaldo19}, and Orthocenter \cite{ronaldo19} were all ellipses. Such a stream of results ushered us, some eight years later, into this second exploratory cycle. Never could we anticipate the surprises still in store!

\begin{figure}
    \centering
    \includegraphics[width=1.0\textwidth]{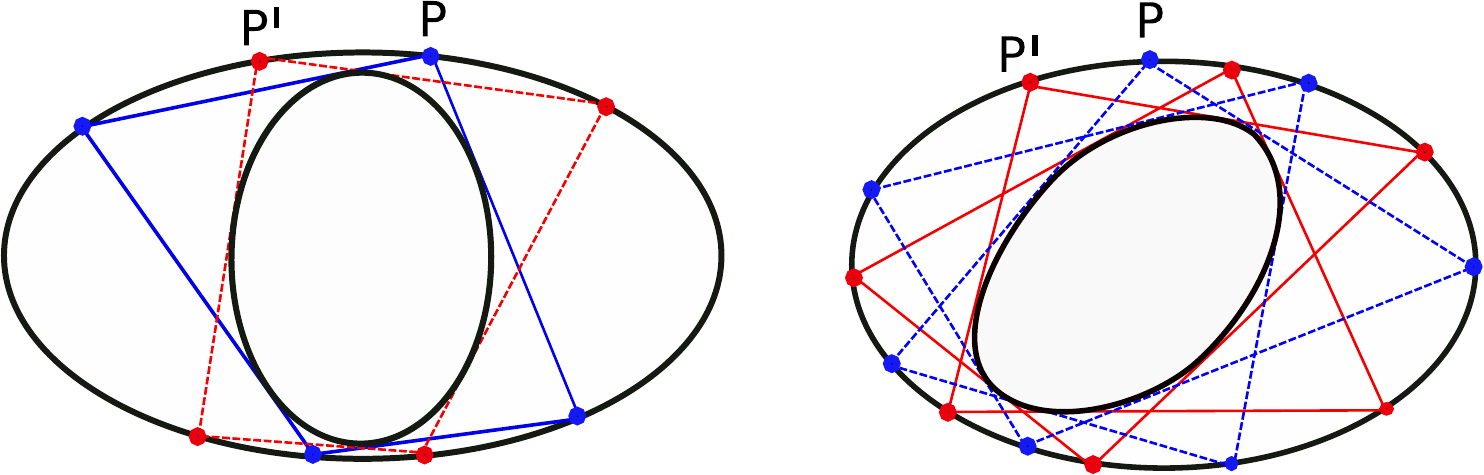}
    \caption{Poncelet's Porism states that given two nested ellipses, if one closed trajectory with $N$ sides can be found starting at a point $P$ on the boundary of the outer
while remaining tangent to the inner, then every boundary point (e.g., $P'$) can initiate an $N$-trajectory, i.e., there exists a one-dimensional {\em family} of closed trajectories.}
\label{fig:poncelet-porism}
\end{figure}

\subsection{An experimental approach and an unexpected breakthrough}

Keen to examine the loci of other {\em Triangle Centers}\footnote{A point defined in terms of vertex positions and angles which is invariant with respect to rigid transformations.} over the $N=3$ orbit family\footnote{Here we also use $N=3$, $N=4$ to denote 3-periodic, 4-periodic, etc.}, we built an interactive \href{https://editor.p5js.org/undefined/present/i1Lin7lt7}{applet} \cite{dsr_applet_x12345}. We started with the first 100 centers catalogued on Clark Kimberling's Encyclopedia\footnote{More than 40,000 are listed there!} \cite{etc}. The loci these produced were quite entertaining: ellipses, circles, a stationary point, a quartic, one with kinks, etc., see our \href{https://dan-reznik.github.io/Elliptical-Billiards-Triangular-Orbits/loci_6tri.html}{locus gallery} \cite{reznik_media}. Interestingly, a few ellipses were similar and/or identical to the Billiard or its Caustic \cite{reznik2020-loci}. 

Beyond loci, the $N=3$ family revealed an amazing property: the ratio of Inradius to Circumradius is 
invariant. Right under our noses, we thought!
In turn, this immediately dictated surprising relations involving orbit angles and areas.

\begin{figure}[ht]
    \centering
    \includegraphics[width=.5\textwidth]{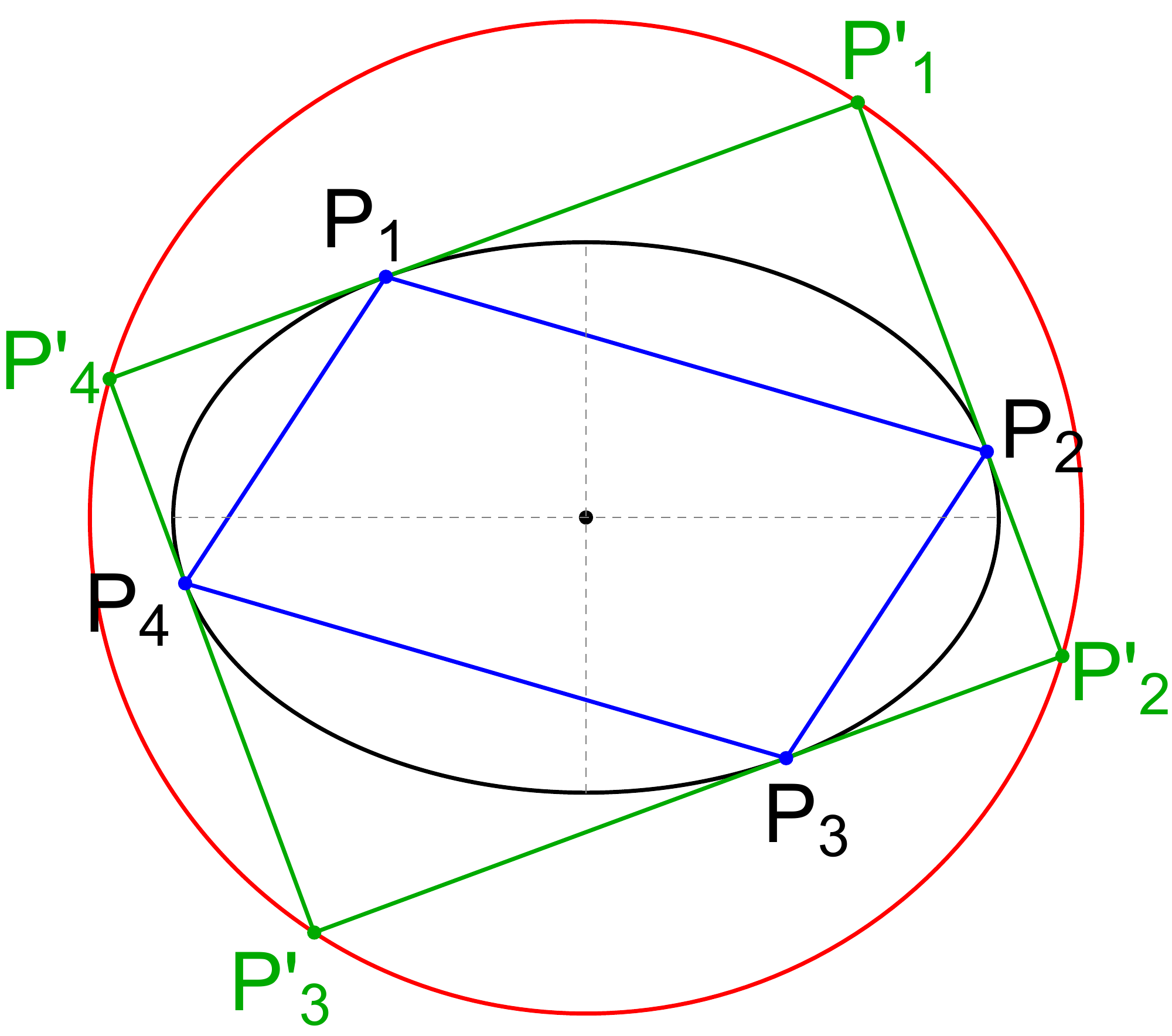}
    \caption{Gaspard Monge (1746--1818) discovered that the locus of points whose tangents to a given ellipse subtend a right angle is a circle (red), called the {\em Orthoptic} Circle. Its connection with Billiards is remarkable \cite{connes07}: From a point $P_1'$ on the circle, shoot two rays tangent to the ellipse, intersecting them with the circle at $P_2'$ and $P_4'$. From either one shoot one new tangent ray and obtain $P_3'$. Four surprises: (i) the intersections form a rectangle (green); (ii) the points of tangency with the ellipse define a parallelogram (blue), which is (iii) a billiard orbit of the ellipse; (iv) all $N=4$ non-intersecting orbits are constant perimeter parallelograms. 
    \href{https://youtu.be/9fI3iM2jrmI}{Video} \cite[pl\#5]{dsr_math_intell_playlist}}    \label{fig:monge-orthoptic}
\end{figure}

At this point a remarkable object steps in, coaxing us into a breakthrough: Monge's Orthoptic Circle,  Figure~\ref{fig:monge-orthoptic}. Its direct connection with the $N=4$ family has been explored by eminent mathematicians \cite{connes07}. Close observation of its geometry provided a bridge with which to generalize $N=3$ invariants to orbits of any $N$. Indeed, there exist generic constructions for point, circular and Caustic loci which are analogous to the $N=3$ case. Additionally, conservation of angular and area quantities already verified in $N=3$ hold true for all $N$. Having started this exploration with the humble triangle, these were remarkable surprises. Merci, Gaspard!

%

As we crossed over to $N>3$, our analytical methods became insufficient. Luckily we were helped by generous mathematicians who contributed insights and proofs \cite{akopyan19_private_meromorphic,helman19,dominique19,olga19_mitten,sergei19_private_circles,sergei19_private_meromorphic}, see Acknowledgements and \cite{sergei2020}. 

We challenge the reader with some questions in the Conclusion and very much welcome feedback. Many of our experiments are available as videos on a YouTube \href{https://bit.ly/2kTvPPr}{playlist} \cite{dsr_math_intell_playlist}. In the text, these are cited as {[n,~pl\#m]}, where n is the standard reference number and m is the video entry into the playlist. Section~\ref{sec:list-videos} provides a quick-reference and links to all videos mentioned below.

\section{Integrability and Conservation}
\label{sec:integrability}
The EB is a 4-dimensional Hamiltonian System, the only planar Billiard known to be integrable \cite{jovanovic11}. This means the particle's path can be explicitly obtained in terms of two Integrals of Motion: (i) Energy (constant velocity with elastic collisions), and (ii) Joachimsthal's (product of angular momenta with respect to the foci \cite{birkhoff1927,sergei91}). The EB can be regarded as a special case of {\em Poncelet's Porism} \cite{dragovic11},  Figure~\ref{fig:poncelet-porism}.

Let the boundary of the EB satisfy $f(x,y)=(x/a)^2+(y/b)^2=1$, where $(a,b)$ are the semi-axes. Continuously, Joachimsthal's Integral implies the trajectory remains tangent to a confocal Caustic (for triangular trajectories, see Figure~\ref{fig:three-orbits}). At the bounces $P_i=(x_i,y_i)$, it implies a quantity $\gamma$ is invariant, given by a remarkable expression:

\begin{equation}
 \gamma=\frac{1}{2}\hat{v}.\nabla{f_i}=\mbox{constant} >0,\,\,\,\mbox{for all }{i}\;,
 \label{eqn:joachim}
\end{equation}

\noindent where $\hat{v}$ is the normalized incoming velocity vector, and:

\begin{equation}
\nabla{f_i}=2\left(\frac{x_i}{a^2}\,,\frac{y_i}{b^2}\right).
\label{eqn:f}
\end{equation}

Constant energy and $\gamma$ constrain the space of trajectories to an abstract 2-torus \cite{sergei91} where the motion has constant frequencies in suitable angular coordinates. When frequencies are rationally-related, the motion is closed and the family has constant period, implying:

\begin{property}
The family of $N$-periodic trajectories has constant perimeter $L$.
\end{property}

Indeed, $L$ and $\gamma$ only depend on $a,b$. For $N=3$ explicit expressions have been derived \cite{ronaldo19,ronaldo19a}.

\section{Locus Pocus: Curves Galore}
\label{sec:loci}
\subsection{Incenter Locus Recap}

\begin{figure}[H]
     \centering
     \begin{subfigure}[t]{0.45\textwidth}
         \centering
         \includegraphics[height=.65 \linewidth]{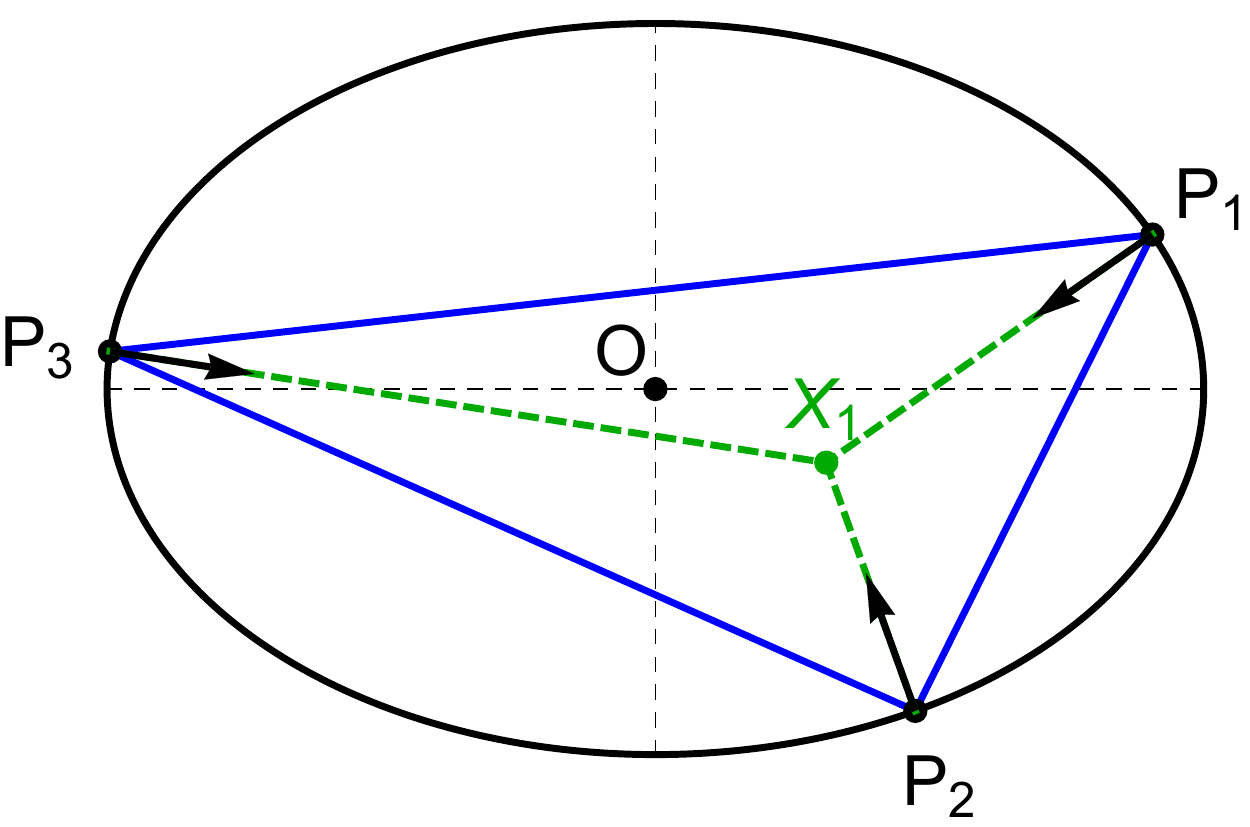}
         \caption{An $N=3$ {\em orbit}, and its Incenter $X_1$: where the bisectors concur.}
         \label{fig:single-orbit}
     \end{subfigure}
     \hfill
     \begin{subfigure}[t]{0.45\textwidth}
         \centering
          \includegraphics[height=.65\linewidth]{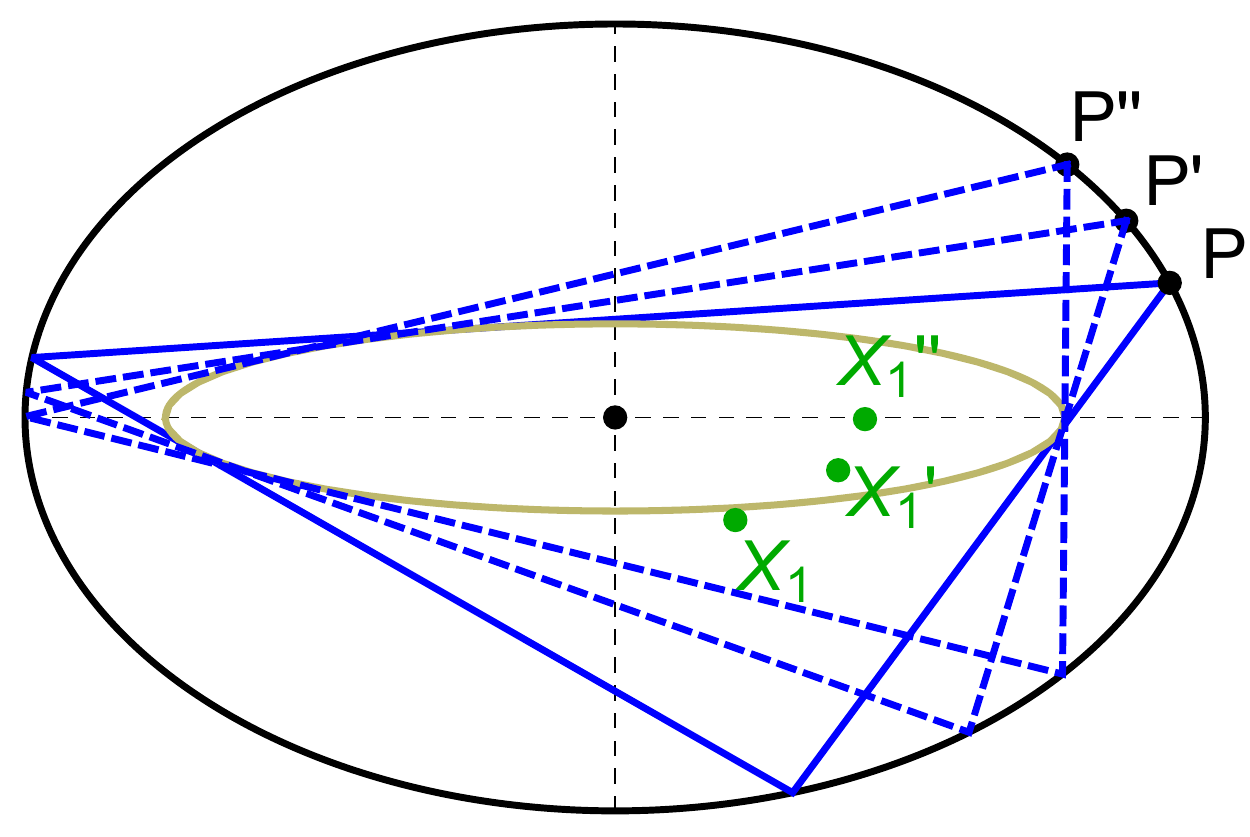}
         \caption{Three triangular orbits, each identified by a starting vertex $P,P',P''$, and their Incenters $X_1,X_1',X_1''$, and the confocal Caustic.
         \label{fig:three-orbits}}
     \end{subfigure}
     \caption{Triangular Orbits in an EB and their Incenter. \href{https://youtu.be/Y3q35DObfZU}{Video} \cite[pl\#6]{dsr_math_intell_playlist}}
      \label{fig:6}
\end{figure}

Consider an EB centered at $O$, and an $N=3$ orbit, shown in Figure~\ref{fig:single-orbit}. A vertex bisector is congruent with the inward normal to the ellipse, and these meet at the orbit's Incenter $X_1$.

Consider the three $N=3$ orbits in Figure~\ref{fig:three-orbits}, identified by a starting vertex $P,P',P''$, as well as their Incenters $X_1,X_1',X_1''$.

The Incenter's Cartesian coordinates in terms of a parametrized vertex, say $P(t)$, are a rather long, non-linear expression \cite{ronaldo19}. This renders remarkable the fact that its locus is an ellipse \cite{olga14} as are those of the Barycenter $X_2$, Circumcenter $X_3$, and Orthocenter $X_4$ \cite{corentin19,ronaldo19,sergei2016}. Recently, the following was also proven \cite{ronaldo19a}:

\begin{theorem}
The locus of the center $X_5$ of the 9-point circle is an ellipse.
\end{theorem}

\begin{figure}[H]
     \centering
     \begin{subfigure}[m]{0.45\textwidth}
     \centering
     \includegraphics[height=.65\linewidth]{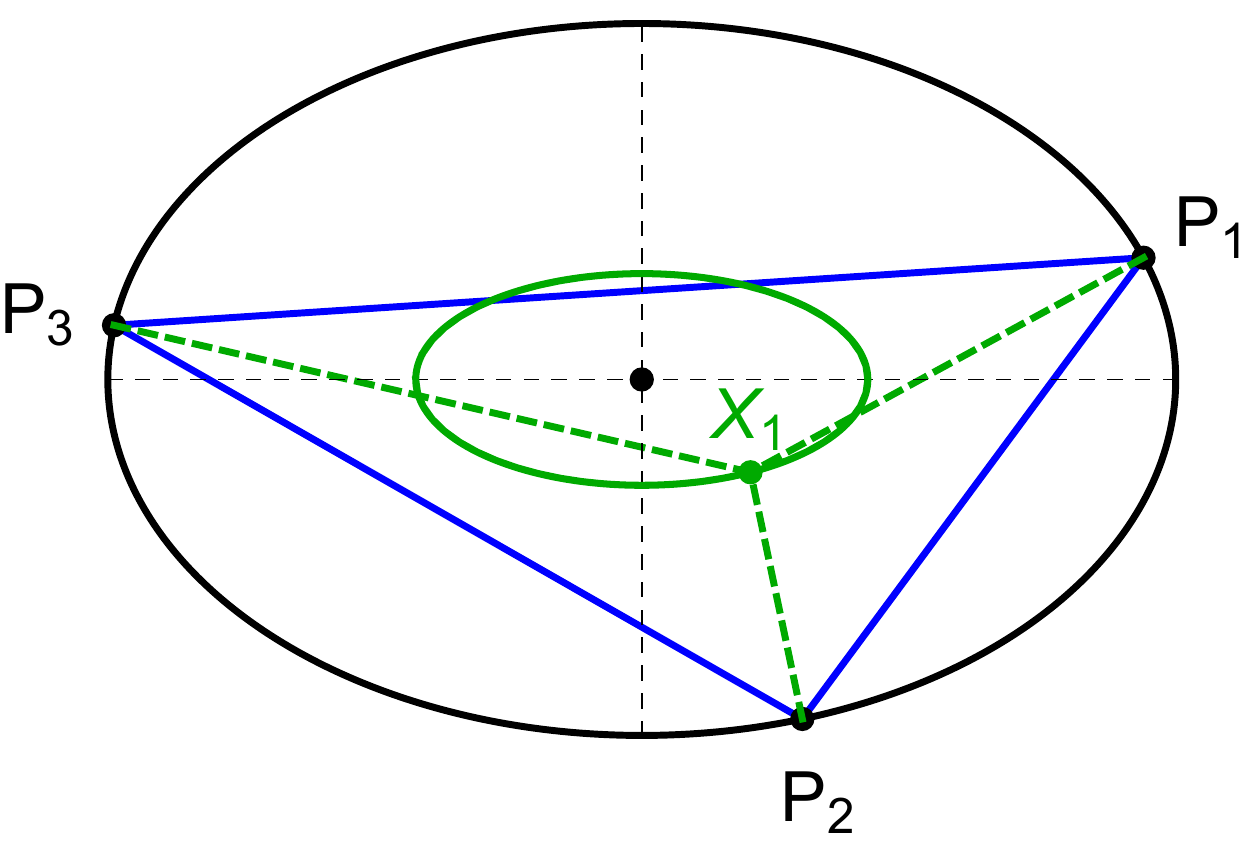}
         \caption{The locus of the Incenter $X_1$ is an ellipse. \href{https://www.youtube.com/watch?v=BBsyM7RnswA}{Video} \cite[pl\#2]{dsr_math_intell_playlist}}
        \label{fig:locus-incenter}
     \end{subfigure}
     \hfill
     \begin{subfigure}[m]{0.45\textwidth}
         \centering
         \includegraphics[height=\linewidth]{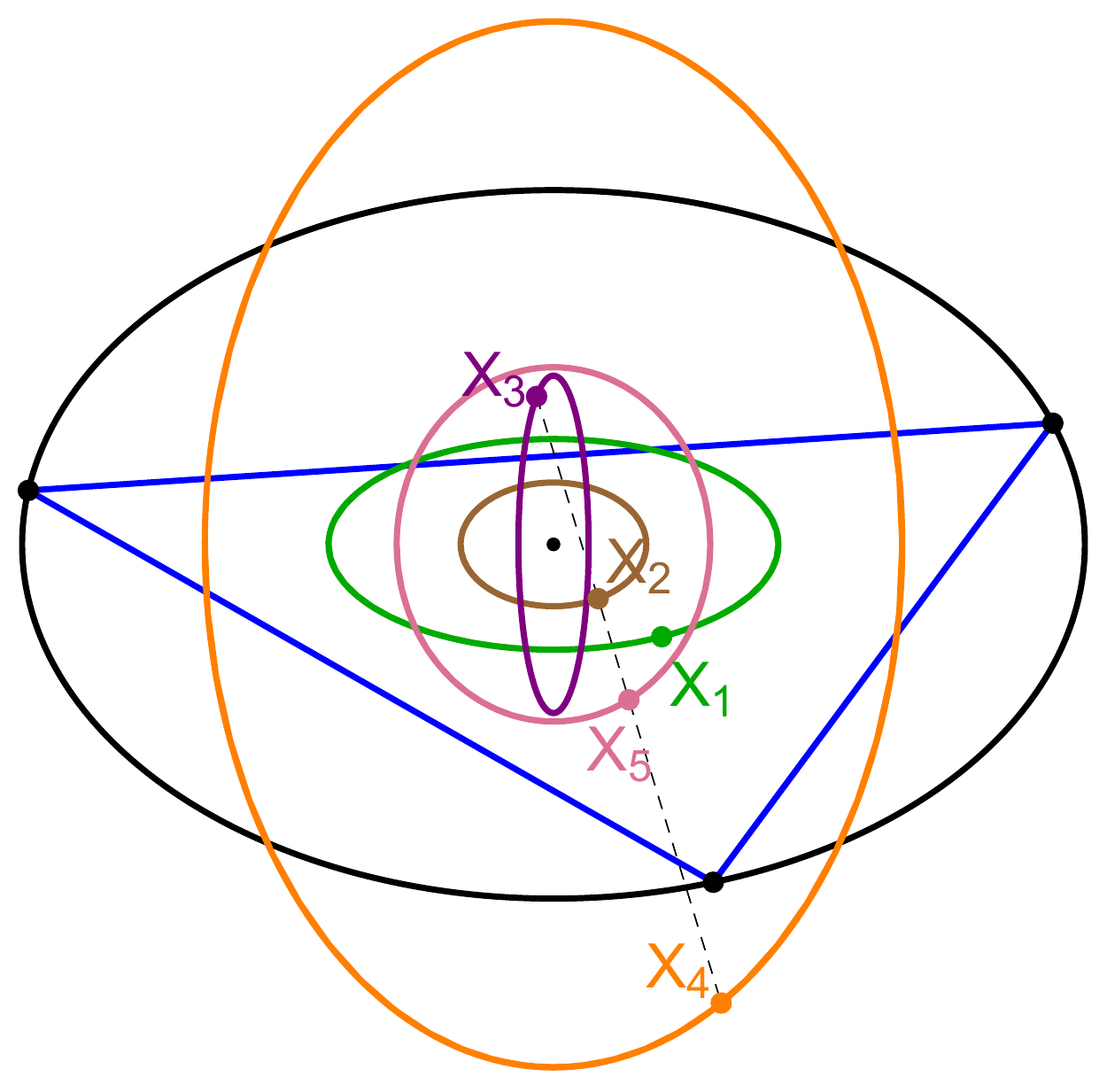}
         \caption{The loci of Incenter $X_1$, Barycenter $X_2$, Circumcenter $X_3$, Orthocenter $X_4$, and Center of the 9-Point Circle $X_5$ are all ellipses. \href{https://youtu.be/sMcNzcYaqtg}{Video} \cite[pl\#7]{dsr_math_intell_playlist}}
         \label{fig:locus-x12345}
     \end{subfigure}
     \caption{Loci of major triangle centers are ellipses.}
     \label{fig:7}
\end{figure}

The above can be seen in Figure~\ref{fig:locus-x12345} and in \cite[pl\#7]{dsr_math_intell_playlist}. Indeed, numerical analysis of the first 100 Kimberling Centers (only 39,900 to go) reported that only 29 of them produce elliptic loci \cite{reznik2020-loci}.

\subsection{Eccentric Excenters}

\begin{figure}[H]
    \centering
    \includegraphics[angle=90,width=.66\textwidth]{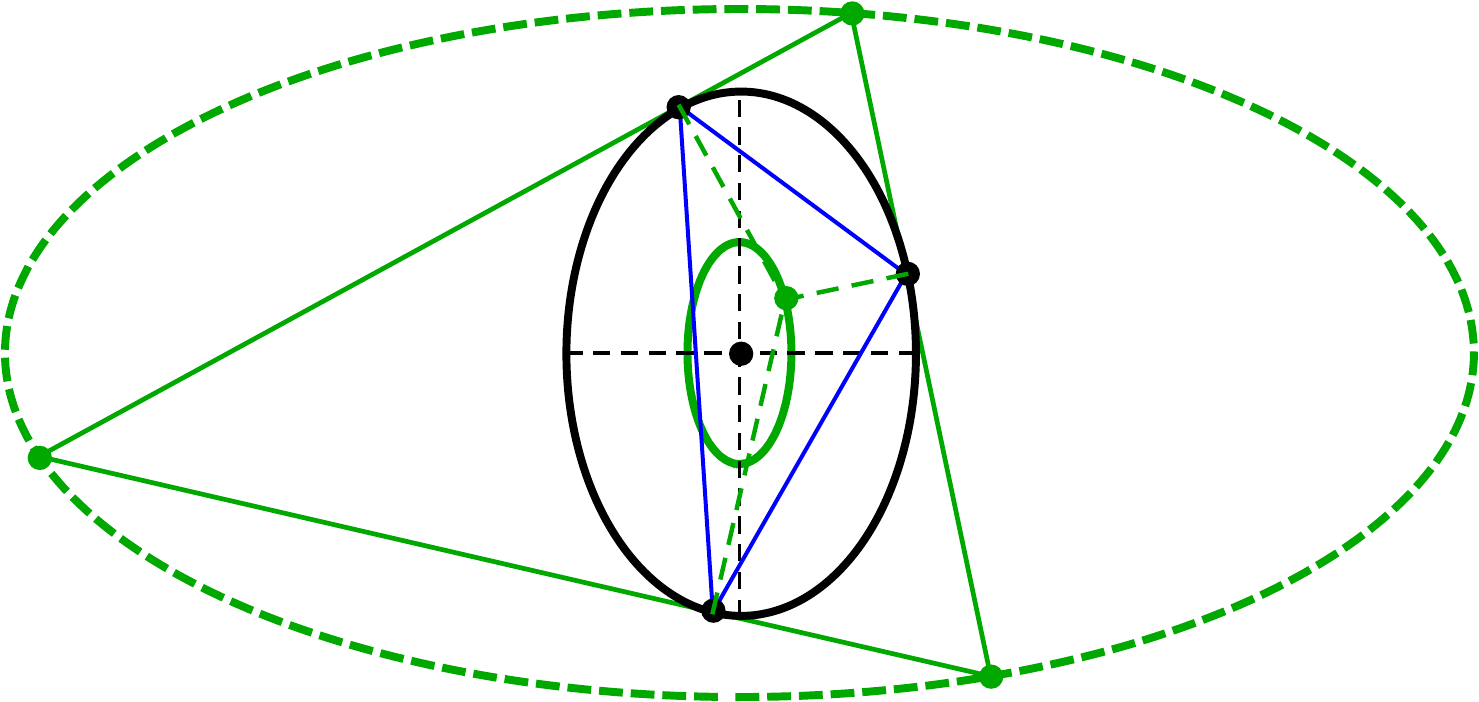}
    \caption{An EB (black) is shown with its axes rotated (to save space), as well as an $N=3$ orbit (blue) and its Excentral Triangle (solid green). The locus of Excenters (vertices of the Excentral Triangle) is an ellipse similar to a perpendicular copy of the locus of the Incenter (shown solid green inside the EB). \href{https://youtu.be/Xxr1DUo19_w}{Video} \cite[pl\#8]{dsr_math_intell_playlist}}
    \label{fig:locus-incenter-excenter}
\end{figure}

The {\em Excenters} are defined in Figure~\ref{fig:constructions}. As depicted in Figure~\ref{fig:locus-incenter-excenter}, it can be shown that \cite{ronaldo19}:

\begin{theorem}
The locus of the Excenters is an ellipse similar to a rotated version of the Incenter's.
\end{theorem}

\begin{figure}[H]
    \centering
    \includegraphics[height=.4\linewidth]{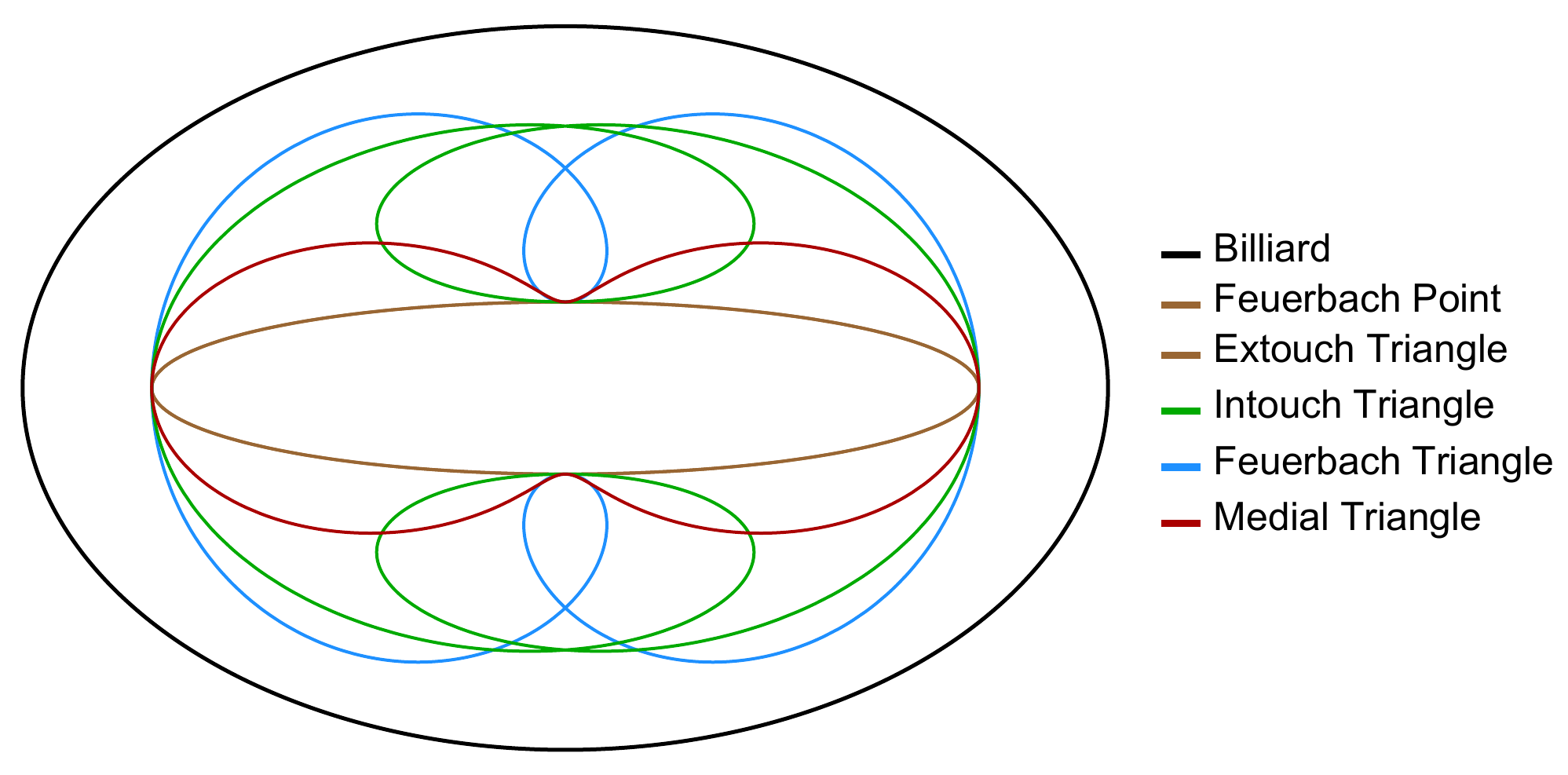}
    \caption{The vertices of the Intouch (green), Feuerbach (blue), and Medial (red) Triangles produce non-elliptic loci. Surprisingly, the Extouchpoints as well as the Feuerbach Point $X_{11}$ (both shown brown), are identical to the $N=3$ Caustic.
    \href{https://youtu.be/OGvCQbYqJyI}{Video} \cite[pl\#9]{dsr_math_intell_playlist}}
    \label{fig:non-elliptic}
\end{figure}

The Excentral Triangle is an example of a {\em Derived Triangle}, i.e., its vertices are computed taking the orbit as a reference triangle. Generally, we've found that vertices of such triangles produce non-elliptic loci, with the Excentral Triangle being an exception. Consider the locus of the {\em Intouch Triangle}.
As mentioned above, their locus is a self-intersecting sextic \cite[pl\#2]{dsr_math_intell_playlist}. In the same vein, vertices of the {\em Feuerbach Triangle}
and {\em Medial Triangle}
both produce non-elliptic loci, Figure~\ref{fig:non-elliptic}. Surprisingly \cite{reznik2020-loci}:

\begin{theorem}
The locus of the vertices of the {\em Extouch Triangle}, where Excircles touch the orbit sides, is an ellipse identical to the Caustic.
\label{thm:extouch}
\end{theorem}

\subsection{Fiery Feuerbach}

\begin{figure}[H]
    \centering
    \includegraphics[width=.60\textwidth]{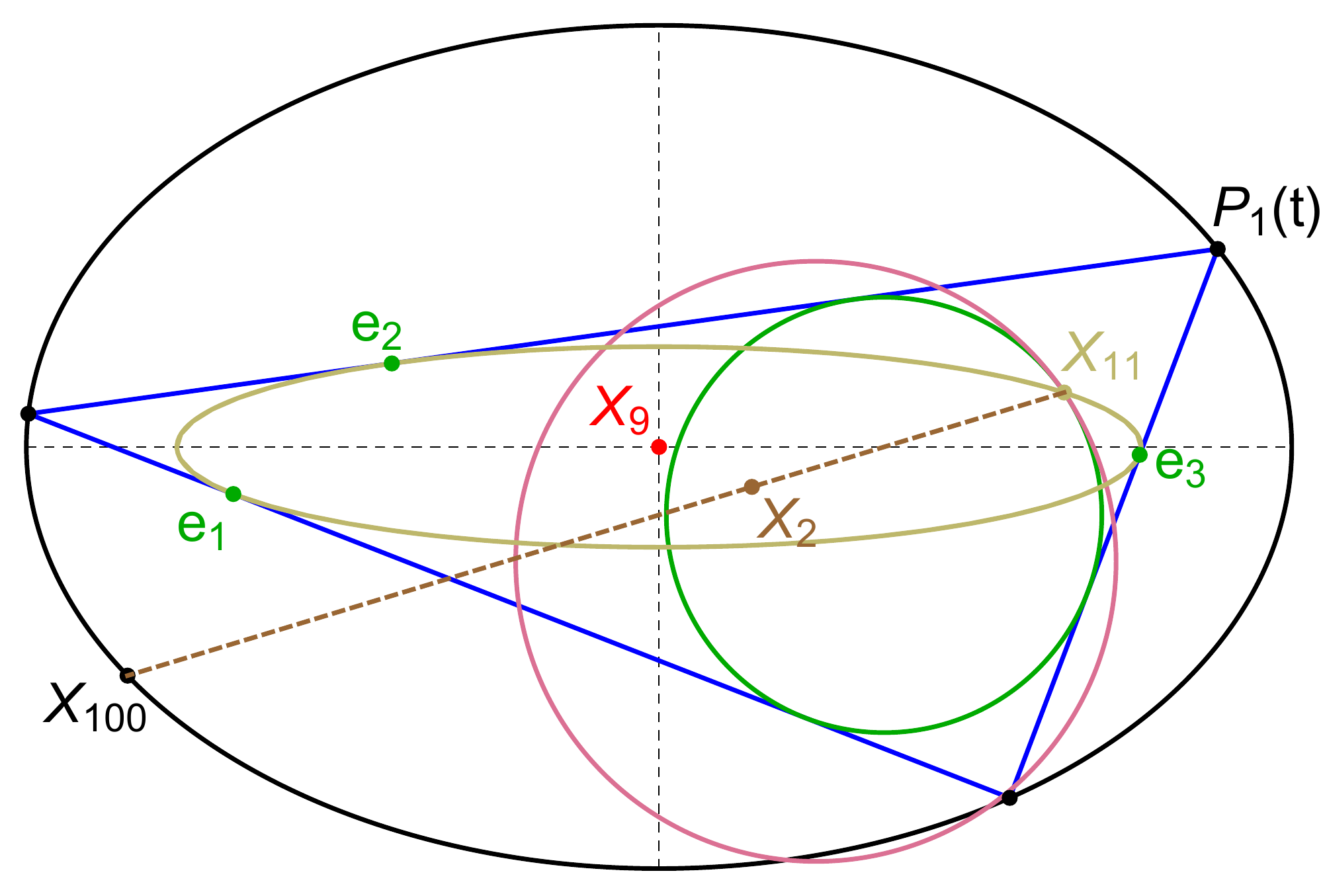}
    \caption{The EB (black), as well as an $N=3$ orbit (blue), with $P_1$ the starting vertex, and the Caustic (brown). Also shown are the orbit's Incircle (green) and 9-Point Circle (pink), whose single point of contact is the Feuerbach Point $X_{11}$. Shown also is its anticomplement $X_{100}$, and the three Extouchpoints $e_1,e_2,e_3$. Remarkable properties include: $X_{11}$ and the Extouchpoints sweep the Caustic, and $X_{100}$ sweeps the Billiard (in opposite directions).
    \href{https://youtu.be/TXdg7tUl8lc}{Video} \cite[pl\#12]{dsr_math_intell_playlist}}
    \label{fig:feuer_loci}
\end{figure}

The Feuerbach Point $X_{11}$  is shown in Figure~\ref{fig:constructions}. $X_{100}$ is its {\em anticomplement}\footnote{A point's double-length reflection about the Barycenter $X_2$.}. As depicted in Figure~\ref{fig:feuer_loci}, 
this duo produces 
a striking phenomenon \cite{reznik2020-loci}:

\begin{theorem}
The locus of the Feuerbach Point $X_{11}$ is identical to the $N=3$ Caustic and the locus of its anticomplement $X_{100}$ is identical to the Billiard.
\end{theorem}

We also noticed that if orbit vertices slide along the Billiard in one direction, $X_{11}$ (resp. the Extouchpoints) will move along the Caustic in the opposite (resp. same) direction, depicted on this \href{https://youtu.be/TXdg7tUl8lc}{video} \cite[pl\#10]{dsr_math_intell_playlist}.

\noindent Other interesting loci worthy of mention include:

\begin{itemize}
    \item Symmedian Point\footnote{Point of concurrence of a triangle's {\em symmedians}, i.e., the reflection of medians about the bisectors.} $X_6$: a convex quartic. When $1<a/b<2$ it closely approximates a perfect ellipse.
    \item Orthic's Incenter: a piecewise-elliptic locus with 4 kinks, shown in a  \href{https://youtu.be/3qJnwpFkUFQ}{video} \cite[pl\#11]{dsr_math_intell_playlist}
    \item Intouchpoints of the Anticomplementary Triangle: identical to the Billiard \cite{minevich17}, shown in a \href{https://youtu.be/50dyxWJhfN4}{video} \cite[pl\#12]{dsr_math_intell_playlist} 
\end{itemize}

\noindent The reader can easily observe the above with our \href{https://editor.p5js.org/undefined/present/i1Lin7lt7}{applet} \cite{dsr_applet_x12345}. 

\section{A Point Locus}
\label{sec:mitten}
\begin{figure}[H]
     \centering
    \includegraphics[width=.8\linewidth]{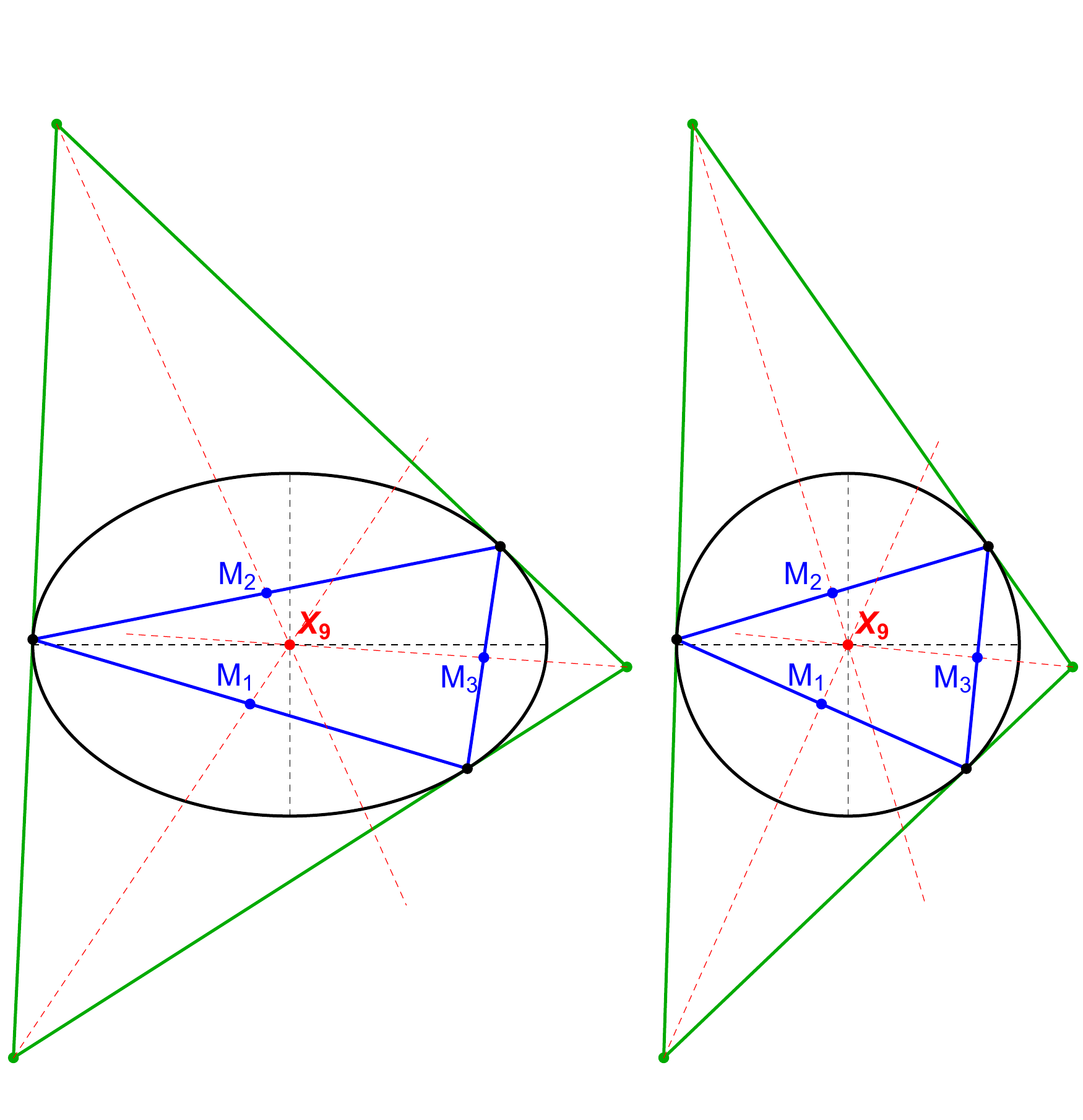}    
     \caption{Affine proof (pun intended!) \cite{olga19_mitten} for the Stationarity of the Mittenpunkt. \textbf{Left}: $N=3$ orbit (blue), its Excentral Triangle (green), and the construction for $X_9$, the Mittenpunkt, stationary at the center of the Billiard (black). \textbf{Right}: Under a suitable affine transform, the Billiard becomes a circle, side midpoints become {\em chord} midpoints (distance ratios are preserved). By symmetry, lines connecting transformed Excenters to the latter must concur at the circle center, uniquely identified under the transform with the Billiard center.
     \href{https://youtu.be/tMrBqfRBYik}{Video} \cite[pl\#13]{dsr_math_intell_playlist}}
     \label{fig:mitten} 
\end{figure}

Triangular centers and vertices of derived triangle produce, we've seen, beautiful loci. But one center, the {\em Mittenpunkt} $X_9$, took us aback. As shown is Figure~\ref{fig:mitten} and on this \href{https://youtu.be/tMrBqfRBYik}{video} \cite[video \#13]{dsr_math_intell_playlist}:

\begin{theorem}
For the 3-periodic family , the Mittenpunkt $X_9$ is stationary at the Billiard's center.
\label{thm:olga}
\end{theorem}

This fact can be proven by an affine transformation \cite{olga19_mitten} explained in Figure~\ref{fig:mitten}. 

A triangle's {\em Circumellipse} \cite{mw} passes through its three vertices. For $N=3$ orbits, we call this object a {\em Circumbilliard}. Since the orbit's Mittenpunkt $X_9$ is congruent with the Billiard's center, we can state (\href{https://youtu.be/vSCnorIJ2X8}{Video} \cite[pl\#14]{dsr_math_intell_playlist}):
\begin{corollary}
Every triangle has a circumellipse to which it is a Billiard orbit, whose center coincides with the triangle's Mittenpunkt.
\end{corollary}

\section{Constant Ratio of Radii}
\label{sec:cosines}
In this Section we describe how we've stumbled upon an incredible manifestation of perimeter and $\gamma$ constancy, involving simple measures in a triangle, and its surprising corollaries.

\begin{figure}[H]
    \centering
    \includegraphics[width=.65\textwidth]{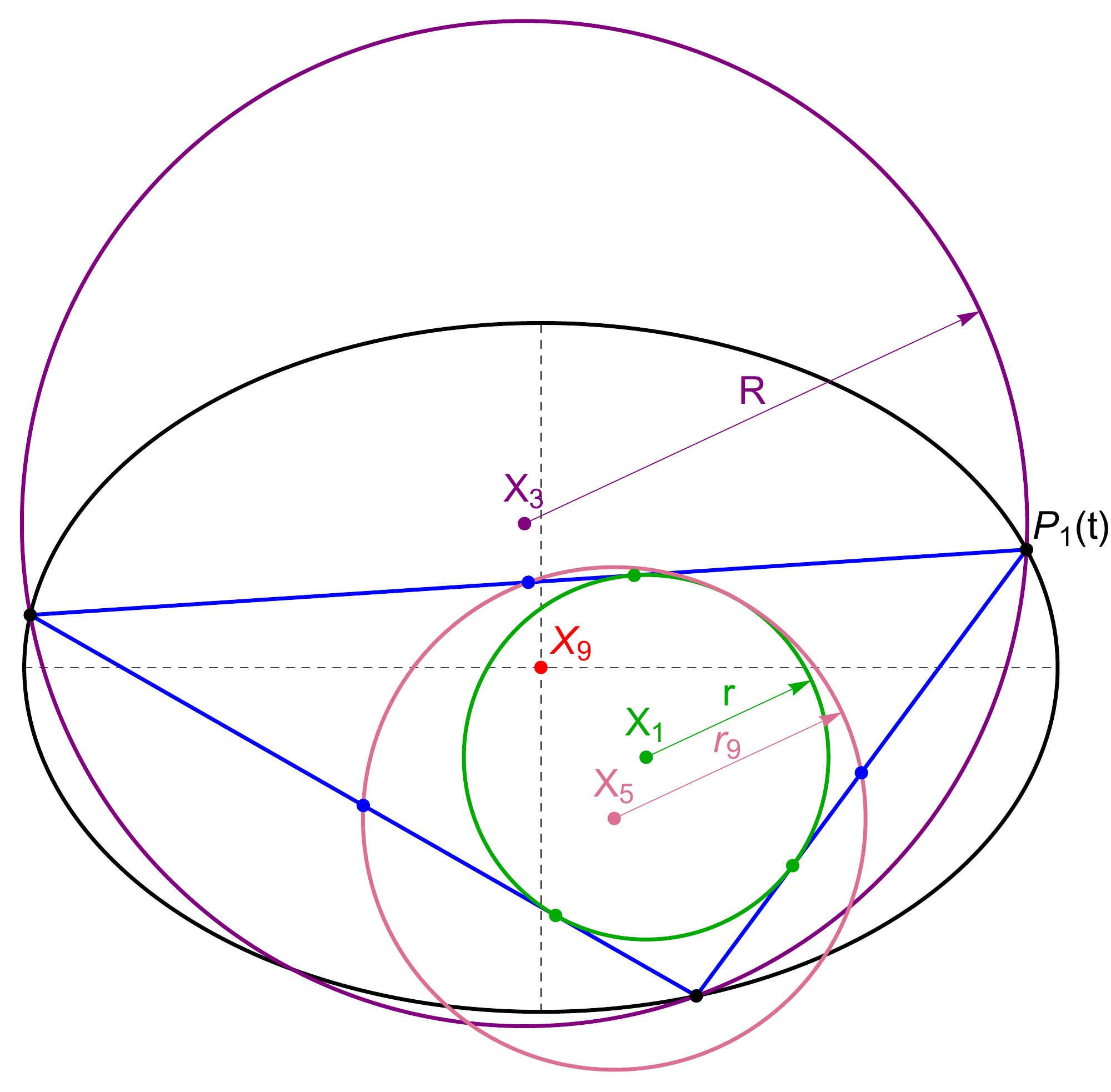}
    \caption{A 3-periodic (blue), a starting vertex $P_1$, the Incircle (green), Circumcircle (purple), and 9-Point Circle (pink), whose centers are $X_1$, $X_3$, and $X_5$, and radii are the inradius $r$, circumradius $R$, and 9-Point Circle radius $r_9$.}
    \label{fig:radii}
\end{figure}

Denote the radii of an $N=3$ orbit's Incircle, Circumcircle, and 9-Point Circle \cite{mw} as $r$ (Inradius) $R$ (Circumradius) and $r_9$, respectively,  Figure~\ref{fig:radii}.

\noindent The relation $R/r_9=2$ is well-known \cite{mw}. Remarkably \cite{ronaldo19a}:

\begin{theorem}
A 3-periodic family conserves the $r/R$ ratio. 
\label{thm:rovR}
\end{theorem}

Specifically, $r/R$ can be expressed in terms of two known constants of motion: perimeter and angular momentum \cite{dominique19,sergei19_private_circles}, as follows:

\begin{equation}
\frac{r}{R}=
\gamma L - 4\, .
\label{eqn:rR_dominique}
\end{equation}

Additionally, it can be shown that when the EB is circular, $r/R$ is maximal at $0.5$ (orbits are equilaterals). As $a/b\rightarrow\infty$, the ratio goes to zero.

\subsection{Conservation of the Sum of Cosines}

Denote $\theta_i,i=1,2,3$ the angles internal to the orbit. The following is an identity valid for any triangle \cite{johnson29}:
\begin{equation}
\sum_{i=1}^{3}{\cos\theta_i}=1+\frac{r}{R}\;.
\label{eqn:rR_cos}
\end{equation}

Since the right-hand side is constant, so must be the sum of the cosines! So a corollary to Theorem~\ref{thm:rovR} is:

\begin{corollary}
\label{cor8}
A 3-periodic family conserves the sum of cosines. 
\end{corollary}

\subsection{Conservation of the Product of Excentral Cosines}

A triangle is always the Orthic of its Excentral \cite{mw}. If $\theta_i'$ are the Excentral's angles, the following is a known relation \cite{johnson29}:

\begin{equation}
\prod_{i=1}^{3}{|\cos\theta_i'|}=\frac{r}{4R}\;.
\label{eqn:prod-cos}
\end{equation}

Let $\theta_i'$ be an angle of the Excentral Triangle opposite to orbit angle $\theta_i$. It can be shown that $\theta_i'=\frac{\pi-\theta_i}{2}$, i.e., the Excentral Triangle is acute. Therefore the absolute sign in Equation~\ref{eqn:prod-cos} can be dropped. Since the right hand side of the above is constant:

\begin{corollary}
\label{cor9}
A 3-periodic family conserves the product of Excentral cosines. 
\end{corollary}

\begin{figure}
    \centering
    \includegraphics[width=\textwidth]{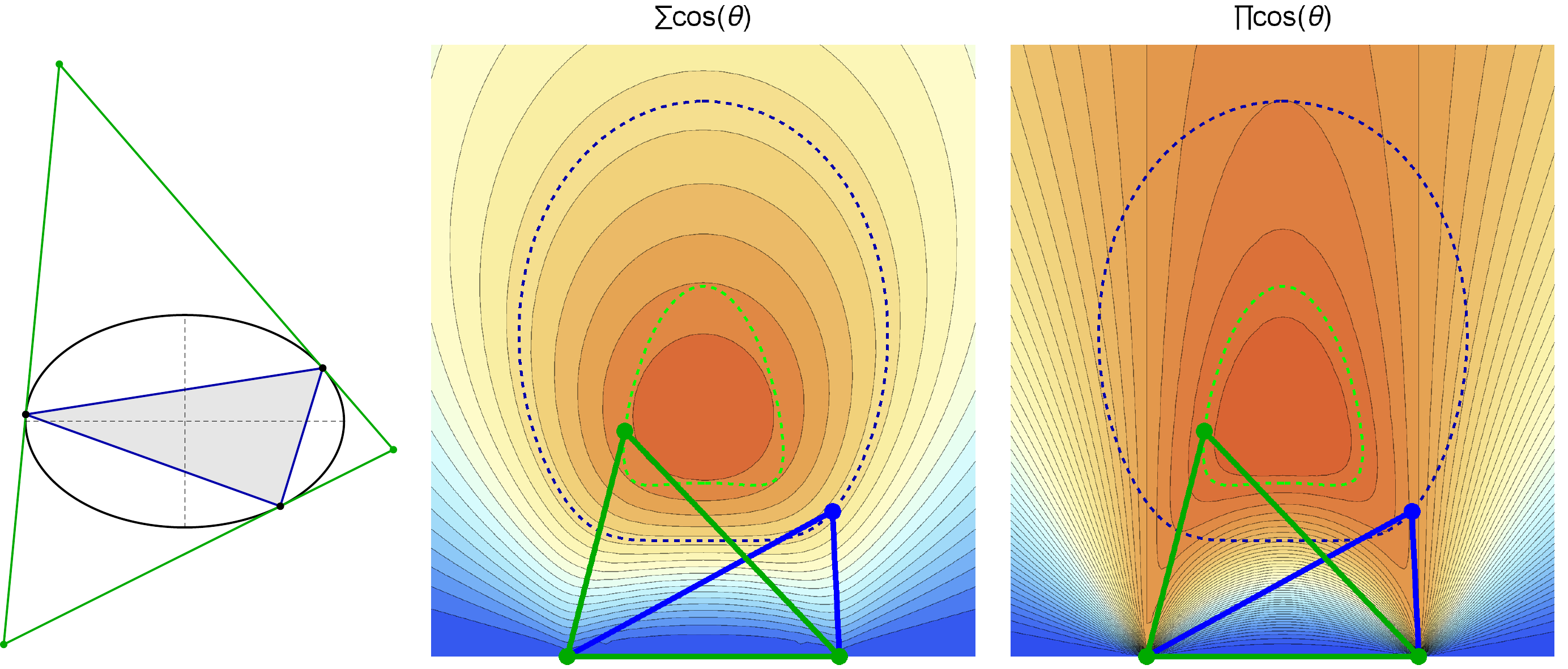}
    \caption{\textbf{Left}: the EB (black), an $N=3$ orbit (blue), and its Excentral Triangle (green). \textbf{Middle}: both the orbit and its Excentral are shown in an ``elementary'' configuration: two vertices are pinned to the $x$ axis and a third one is set free, in such a way that triangles here are similar to their counterparts on the left. Isocurves of constant cosine sum are shown in terms of the position of the free vertex. Notice only the orbit triangle follows such an isocurve (dashed blue). \textbf{Right}: Here isocurves of cosine product are shown, showing the free vertex of the Excentral Triangle moves along one (dashed green).  \href{https://youtu.be/P8ykpE_ZbZ8}{Video}  \cite[pl\#15]{dsr_math_intell_playlist}}
    \label{fig:conserve_cosines}
\end{figure}

\noindent Figure~\ref{fig:conserve_cosines} illustrates conservation of sum (resp. product) of $N=3$ orbit (resp. Excentral Triangle) cosines.

\subsection{Conservation of Excentral-to-Orbit Area Ratio}

Let $A$ be the area of some triangle and $A_h$ the area of its Orthic. A well-known relation is that the ratio $A/A_h$ is inversely proportional to the Inradius-to-Circumradius ratio of the Orthic \cite{johnson29}:

\begin{equation}\label{eqn:AAh}
\frac{A}{A_h}=\frac{2R_h}{r_h}\;. 
\end{equation}

\noindent Since the orbit is its Excentral's Orthic, $r/R$ constant implies:

\begin{corollary}
\label{cor10}
A 3-periodic family conserves the Excentral-to-Orbit area ratio.
\end{corollary}

\section{A Circular Locus}
\label{sec:circles}
Loci can be ellipses, but can they be a circle? Let $P$ be an $N=3$ vertex and $P'$ its reflection about the origin. By symmetry, this will ``fall" on the Billiard. Let  $Q_1$ and $Q_2$ be intersections of the tangent to the Billiard at $P'$ with the orbit's Excentral Triangle $T'$. The following remarkable property holds, illustrated in Figure~\ref{fig:cosine_circle_locus}: 

\begin{figure}[H]
\centering
\includegraphics[width=.8\linewidth]{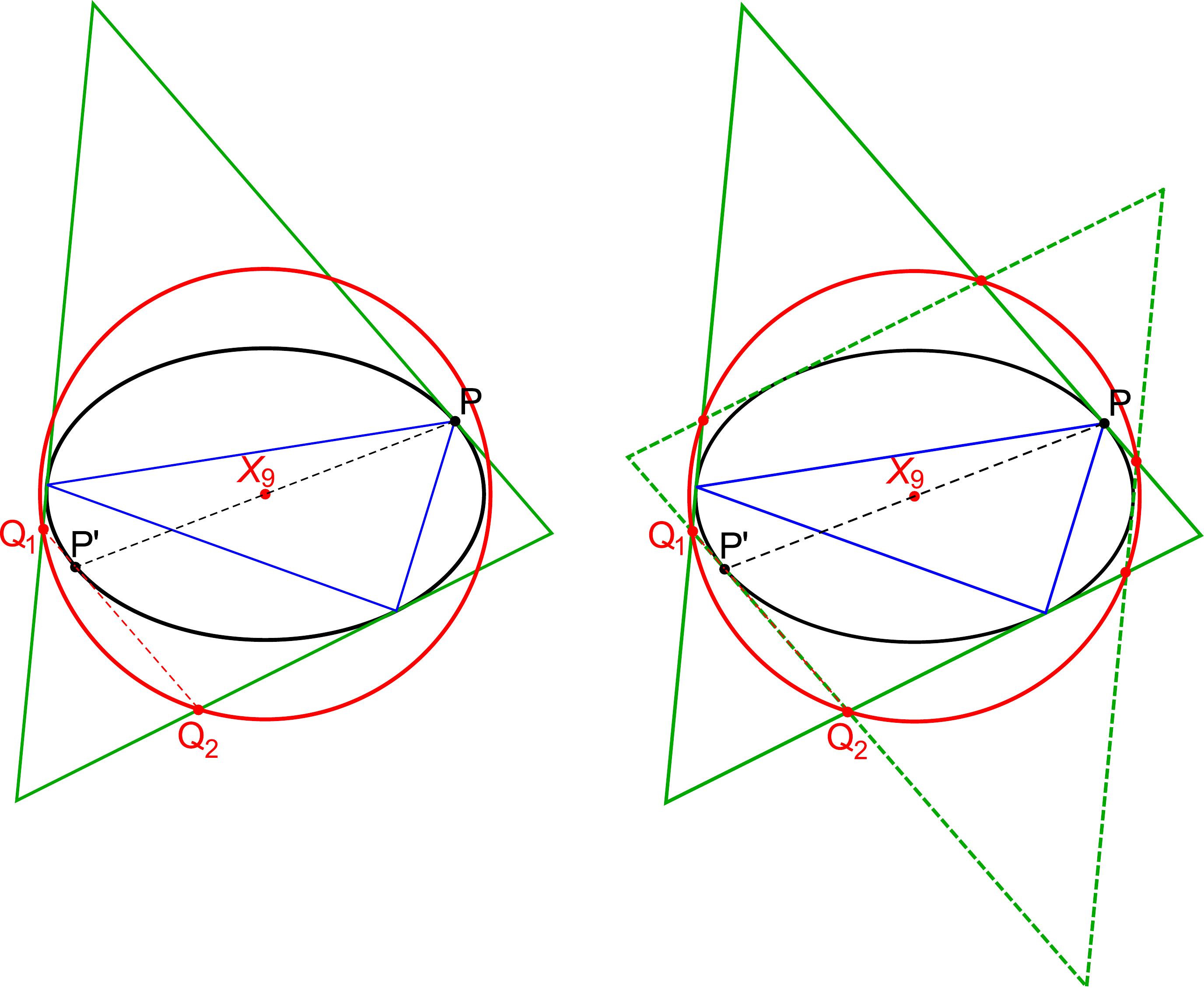}
\caption{\textbf{Left}: Circular Locus of $Q_1$ (or $Q_2$). \textbf{Right}: The Excentral Triangle's six intersections with its reflection about its Symmedian Point (congruent with $X_9$) is a circle. \href{https://youtu.be/CrOSI8d8qDc}{Video 1}, \href{https://youtu.be/hCQIT6_XhaQ}{Video 2} \cite[pl\#16,17]{dsr_math_intell_playlist}}
\label{fig:cosine_circle_locus}
\end{figure}

\begin{theorem}
The locus of $Q_1$ (or $Q_2$) is a circle centered on the Billiard's.
\end{theorem}

The circular locus\footnote{For Triangle enthusiasts, the circular locus is congruent with the {\em Cosine Circle} of the Excentral Triangle, also known as the {\em Second Lemoine Circle}, a kind of {\em Tucker Circle} \cite{mw}. Its center is the Symmedian Point of the Excentral. Since the latter is congruent with the orbit's Mittenpunkt, both radius and center are stationary. No other Tucker Circles for the $N=3$ family have been found to be stationary.} is always external to the Billiard and its radius $r^*$ can be written as a function of the aspect ratio \cite{ronaldo19a}, or even more simply \cite{dominique19,sergei19_private_circles}:

\begin{equation}
r^* = 
\frac{L}{\frac{r}{R}+4} = \frac{1}{\gamma}\;.
\label{eqn:rstar}
\end{equation}

This surprising phenomenon can be viewed here \cite[pl\#16,17]{dsr_math_intell_playlist}.

\section{Monge Madness: onward to N{\textgreater}3}
\label{sec:generalize}
Close observation of the $N=4$ orbit family, Figure~\ref{fig:monge-orthoptic}, provided invaluable clues with which to generalize $N=3$ invariants to $N>3$ orbits. Namely:

\begin{itemize}
    \item Mittenpunkt: The lines connecting the Tangential Polygon's vertices to the parallelogram midpoints concur at the Ellipse's center.
     \item Null Sum of Cosines: The polygon formed by the four 
     tangency points is a parallelogram therefore the sum of its cosines is zero.
    \item Null Product of Cosines: The Tangential Polygon is a rectangle and therefore the product of its cosines is zero.
   \item Circular Locus: The orthoptic locus is a circle.
\end{itemize}

\subsection{$N>3$ Stationary Point}

Let the Tangential Polygon be the generalization of the Excentral Triangle for $N>3$. The same Affine proof used in Theorem~\ref{thm:olga} and Figure~\ref{fig:mitten} yields:

\begin{theorem}
Given an $N\geq{3}$ orbit family, the point of concurrence of lines drawn from vertices of the tangential polygon through the midpoints of orbit sides is stationary at the center of the Billiard.
\end{theorem}

\begin{figure}
    \centering
    \includegraphics[width=\textwidth]{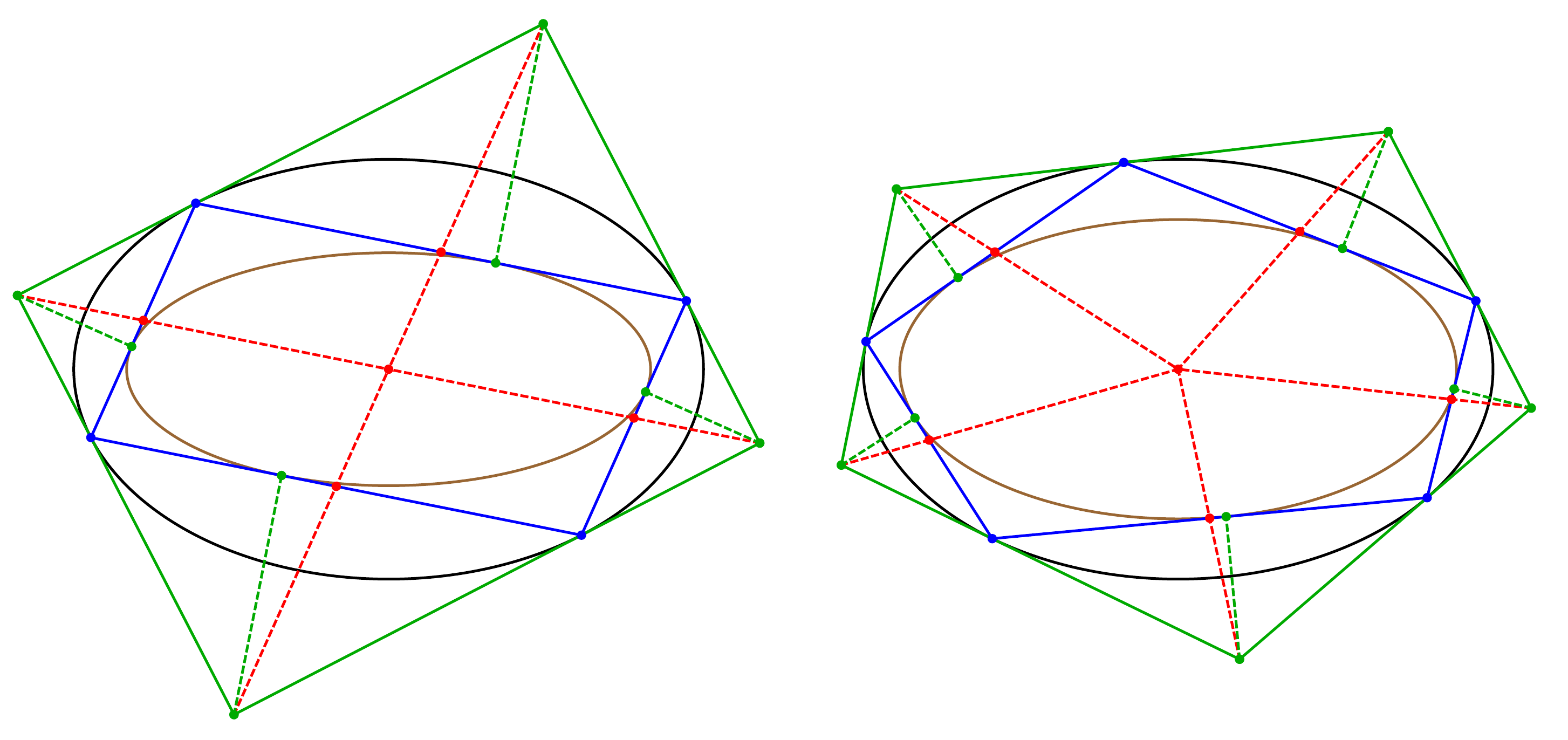}
    \caption{\textbf{Left}: $N=4$ (resp. \textbf{Right}: $N=5$) orbits shown in blue. The dashed red lines connect vertices of the Tangential (green) to orbit sides' midpoints. Notice they intersect at the center of the Billiard, a type of generalized Mittenpunkt. Also shown are (dashed green) perpendiculars dropped from Tangential Vertices onto orbit sides. Their feet lie on and sweep the Caustic (brown). \href{https://youtu.be/TV2p7fPlYfE}{Video 1}, \href{https://youtu.be/Bpc-MrR2IMc}{Video 2}  \cite[pl\#18,19]{dsr_math_intell_playlist}.}
    \label{fig:gen-mitten}
\end{figure}

\noindent This is illustrated in Figure~\ref{fig:gen-mitten} and in   \cite[pl\#13]{dsr_math_intell_playlist}.

\subsection{Generalized Extouchpoints}

Theorem~\ref{thm:extouch} states Extouchpoints are on and sweep the $N=3$ caustic. Analogously, Figure~\ref{fig:gen-mitten}:

\begin{theorem}
The feet of perpendiculars dropped from vertices of the Tangential Polygon onto orbit sides are on the caustic and have the latter as their locus.
\end{theorem}

It turns out the above is a special case of Chasles' Theorem \cite{sergei2016proj}. 

This offers a method to easily ``find'' a point on the caustic which only requires two orbit vertices, alternate to \cite{himmelstrand12} which required the Billiard's foci.

\subsection{$N>3$ Cosine Sum and Product}

The fact that both $N=3$ and $N=4$ conserve cosine sum suggests $N>4$ might as well, and this was first confirmed numerically for $N=5,\ldots,30$ non-intersecting orbits at various Billiard aspect ratios. A proof to this surprising fact has been kindly contributed \cite{sergei19_private_meromorphic}.

\begin{theorem}
The sum of cosines is conserved for non-intersecting orbits, for all $N$. 
\label{thm:cosine-sum-meromorphic}
\end{theorem}

\noindent Beyond establishing the sum of cosines was constant, we noticed Equations~\ref{eqn:rR_dominique}~and~\ref{eqn:rR_cos} imply that for $N=3$, $\sum{\cos\theta_i}=\gamma{L}-3$. The appearance of a ``$3$'' suggested we should verify (numerically) if substituting this digit by $N$ would hold for all $N$. Surprisingly this did hold, and a proof soon followed \cite{sergei19_private_cosine_sum_expression}:

\begin{theorem}
For an $N$-periodic orbit, $\sum_{i=1}^{N}{\cos\theta_i}=\gamma{L}-N$.
\label{thm:cosine-sum}
\end{theorem}

\noindent Because all terms $\gamma$, $L$, and $N$ are constant, the above subsumes and is more specific than Theorem~\ref{thm:cosine-sum-meromorphic}. Note that for $N=4$, the sum of cosines is zero and $\gamma{L}=4$, independent of the EB's aspect ratio.

Furthermore, since $\cos\frac{\theta_i}{2}=\frac{\nabla{f_i}.\hat{v}}{||\nabla{f_i}||}=\frac{2\gamma}{||\nabla{f_i}||}$, the following expression for the perimeter $L$ will hold:

\begin{corollary}
\begin{equation}
    L=8\gamma\sum_{i=1}^{N}\frac{1 }{||\nabla{f_i}||^2} \;\cdot
    \label{eqn:perimeter-new}
\end{equation}
\end{corollary}

\noindent Since both $N=3$ and $N=4$ conserve the product of their Excentral/Tangential cosines, we confirmed numerically that $N=5,\ldots,30$ non-intersecting orbits at various aspect ratios also conserve this quantity. It can also be proven that \cite{akopyan19_private_meromorphic}:

\begin{theorem}
The product of cosines for the Tangential Polygon of non-intersecting orbits for {\em any} $N$ is conserved.
\end{theorem}

\subsection{Area Ratio for $N>3$}

For $N=3$ the Excentral-to-Orbit area ratio is conserved, however, we can numerically check this is not the case for $N=4$. Proceeding to $N>4$ we verify that only odd $N$ preserve area ratio, a fact which has been subsequently formally proven \cite{akopyan19_private_meromorphic}.

\begin{theorem}
The ratio of areas between Tangential Polygon and Non-Self-Intersecting Orbit is conserved for all odd $N$.
\end{theorem}

\subsection{$N>4$ Circular Loci}

\begin{figure}[H]
    \centering
        \includegraphics[width=\textwidth]{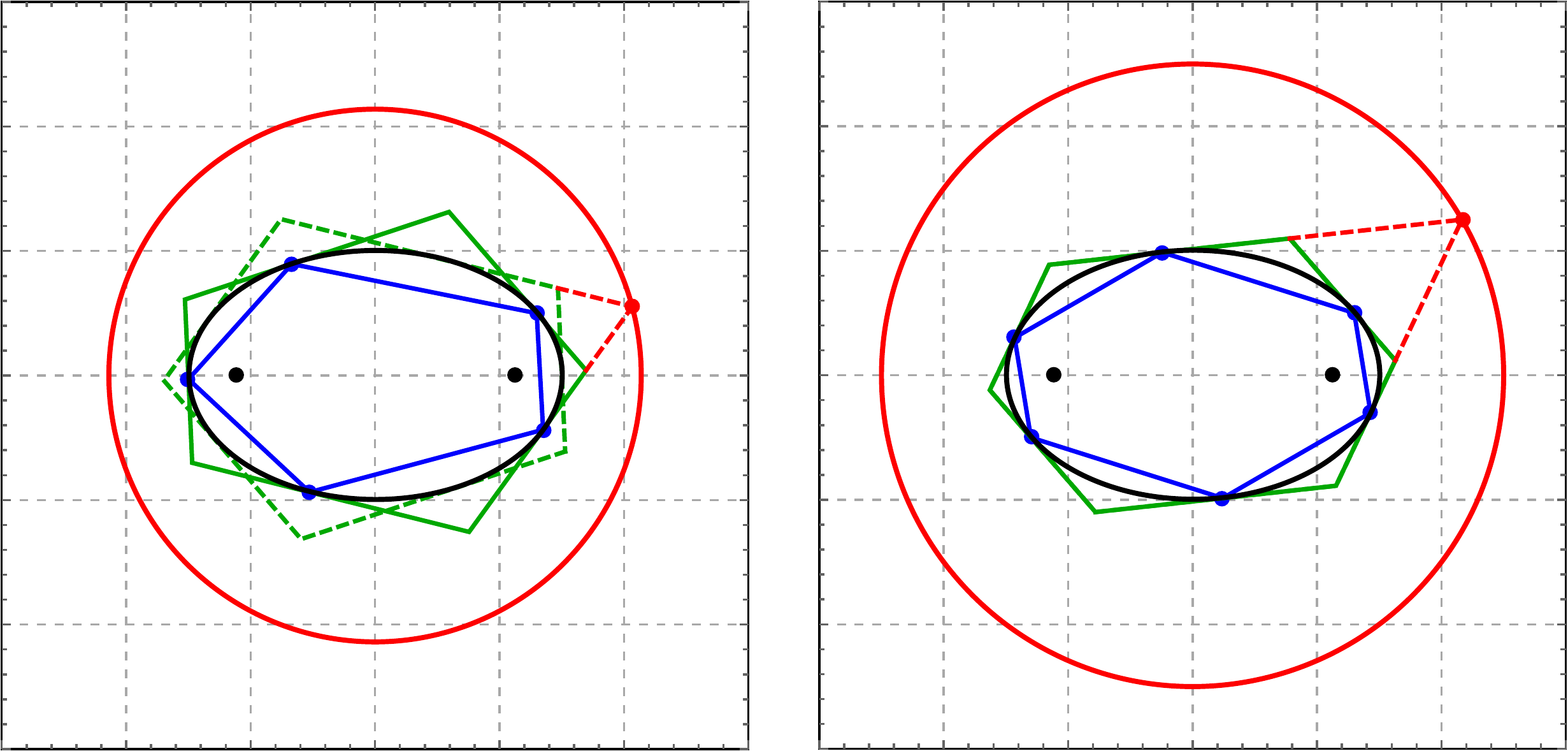}
    \caption{The construction of a circular locus for $N=5,6$.
    \href{https://youtu.be/dINE4aH1cvk}{Video 1},
    \href{https://youtu.be/EFeINGIDFrg}{Video 2} \cite[pl\#20,21]{dsr_math_intell_playlist}}
    \label{fig:gen-circ-grid}
\end{figure}

For the $N=3$ case,  Figure~\ref{fig:cosine_circle_locus}, the locus of the intersection of an Excentral Triangle edge with an alternate edge reflected about the Billiard center is a circle. For $N=4$ we have Monge's Orthoptic Circle. Here is how these two constructions can be unified to $N>4$, Figure~\ref{fig:gen-circ-grid}, \cite[pl\#5]{dsr_math_intell_playlist}:

\begin{theorem}
Let $O$ be the center of the Billiard. The locus of the intersection of an edge of the Tangential Polygon with the reflection of the next tangential edge about $O$ is a circle centered on $O$. Its radius is $r^*=1/\gamma$.
\end{theorem}

\noindent {\bf Proof} \cite{sergei19_private_circles}: Let $\nabla_i$ denote $\nabla{f_i}$ of Equation~\ref{eqn:f}. Consider two consecutive orbit vertices $P_{i}$ and $P_{i+1}$. Momentum conservation (Equation~\ref{eqn:joachim}) implies $\hat{v}.\nabla_i= -\hat{v}.\nabla_{i+1}=2\gamma$. So we have:

\begin{equation}
    \hat{v}.\left(\nabla_i+ \nabla_{i+1}\right)=0\;.
    \label{eqn:sum-nablas}
\end{equation}

\noindent Since $\nabla_i$ (resp. $\nabla_{i+1}$) is normal to the ellipse at $P_{i}$ (resp. $P_{i+1}$), a point $z$ on the tangent line at $P_i$ (resp. $-P_{i+1}$) is given by $z.\nabla_i=2$ (resp. $z.\nabla_{i+1}=-2$). Let $z$ be where both lines intersect, $z.\left(\nabla_i+\nabla_{i+1}\right)=0$. It follows from Equation~\ref{eqn:sum-nablas} that $z$ is parallel to $\hat{v}$. Since $\hat{v}.\nabla_i=2\gamma$ and $z.\nabla_i=2$, we have $z = \hat{v}/\gamma$. That is, $z$ lies on the circle of radius $1/\gamma$.\qed

\section{Conclusion \& Questions}
\label{sec:conclusion}

The reader is invited to check out more details about this work on a companion \href{https://dan-reznik.github.io/Elliptical-Billiards-Triangular-Orbits/}{website} \cite{reznik_web}. We part with the following questions:

\begin{itemize}
    \item For $N=3$, what determines whether a triangle center or derived vertex produces an elliptic vs a more complicated locus? Which loci are algebraic, and which are not? Is there a general theory?
    \item Are there ellipsoidal (3d) counterparts to these invariants?
    \item Which invariants are still true for self-intersecting orbits? (consider videos for \href{https://youtu.be/cCYxN7ueGV4}{N=4} and \href{https://youtu.be/ECe4DptduJY}{N=5} \cite[pl\#3,4]{dsr_math_intell_playlist})
    \item Are there invariants for non-billiard (Poncelet) orbit families, e.g., created with two non-confocal or misaligned ellipse pairs, as shown in \cite{sergei2016}? A recent \href{https://youtu.be/B5dRXT8Xerw}{Video} \cite[pl\#22]{dsr_math_intell_playlist} shows certain aligned Poncelet pairs can collapse the locus of Triangular Centers to a point ($X_1$,$X_2$,$X_3$,$X_4$,$X_6$). Can an explicit map be found where stationary points under some special Poncelet Pair are identified with the Mittenpunkt of a Billiard (confocal) pair?
    \item Are there properties of triangle centers if orbits are defined on the surface of a sphere where edges become arcs of great circles (geodesics)? How about on the surface of an ellipsoid, bounded by the lines of curvature? Such Billiards are also known to be integrable \cite{sergei2002}.
    \item In the spirit of Theorem~\ref{thm:cosine-sum}, can expressions be derived in terms of $\gamma,L,N$ for the constant product of tangential polygon cosines and tangential-to-orbit area ratios? 
\end{itemize}

\subsection{List of Videos}
\label{sec:list-videos}

Videos mentioned above have been placed on a Youtube \href{https://bit.ly/2kTvPPr}{playlist} \cite{dsr_math_intell_playlist}.  Table~\ref{tab:playlist} contains quick-reference links to all videos mentioned, with column ``PL\#'' providing video number within the playlist.

\begin{table}[H]
\begin{tabular}{llll}
Title & Fig. & \href{https://bit.ly/2kTvPPr}{Links} & PL\#\\
\hline
{Open trajectories: varying input angle} &
\ref{fig:billiard-trajectories} &
\href{https://youtu.be/A7mPzrNJHkA}{N=3} &
{1} \\ 

{Locus of Incenter and Intouchpoints} &
\ref{fig:intro-plot},\ref{fig:locus-incenter} &
\href{https://www.youtube.com/watch?v=9xU6T7hQMzs}{N=3} &
{2}\\

{Self-intersecting orbits} &
-- &
\href{https://youtu.be/cCYxN7ueGV4}{N=4}, \href{https://youtu.be/ECe4DptduJY}{N=5} &
{3,4} \\

{Monge's Orthoptic Circle} &
\ref{fig:monge-orthoptic} & 
\href{https://youtu.be/9fI3iM2jrmI}{N=4} &
{5} \\

{Orbits and their Caustics} &
\ref{fig:6} &
\href{https://youtu.be/Y3q35DObfZU}{N=3\ldots{6}}& 
{6} \\

{Elliptic loci for $X_1\ldots{X_5}$} &
\ref{fig:7} & 
\href{https://youtu.be/sMcNzcYaqtg}{N=3} &
{7} \\

{Elliptic locus of the Excenters} &
\ref{fig:locus-incenter-excenter} &
\href{https://youtu.be/Xxr1DUo19_w}{N=3} &
{8} \\

{Loci of Medial, Intouch, Feuerbach Tri's} &
\ref{fig:non-elliptic} &
\href{https://youtu.be/OGvCQbYqJyI}{N=3} &
{9} \\

{Locus of Feuerbach Pt. and Anticompl.} &
\ref{fig:feuer_loci} &
\href{https://youtu.be/TXdg7tUl8lc}{N=3} &
{10} \\
 
{Locus of Orthic Incenter is 4-piece ellipse} &
-- &
\href{https://youtu.be/3qJnwpFkUFQ}{N=3} &
{11} \\

{Anticompl. Intouchpts. are on the Billiard} &
-- &
\href{https://youtu.be/50dyxWJhfN4}{N=3} & 
{12} \\

{Mittenpunkt stationary at Billiard center} & \ref{fig:mitten} &
\href{https://youtu.be/tMrBqfRBYik}{N=3} &
{13} \\

{The Circumbilliard} &
-- &
\href{https://youtu.be/vSCnorIJ2X8}{N=3} &
{14} \\

{Constant cosine sum and product} &
\ref{fig:conserve_cosines} &
\href{https://youtu.be/P8ykpE_ZbZ8}{N=3} &
{15} \\

{Stationary Cosine Circle} &
\ref{fig:cosine_circle_locus} &
\href{https://youtu.be/CrOSI8d8qDc}{v1},{ }
\href{https://youtu.be/hCQIT6_XhaQ}{v2} & 
{16,17} \\

{Generalized Mittenpunkt and Extouchpts.} &
\ref{fig:gen-mitten} &
\href{https://youtu.be/Bpc-MrR2IMc}{N=4,5},{ }
\href{https://youtu.be/TV2p7fPlYfE}{N{\textgreater}3} &
{18,19} \\

{Generalized stationary circle} &
\ref{fig:gen-circ-grid} &
\href{https://youtu.be/dINE4aH1cvk}{N=5},{ } \href{https://youtu.be/EFeINGIDFrg}{N=3\ldots{8}} &
{20,21} \\

{Loci with Poncelet Ellipse Pairs} &
-- &
\href{https://youtu.be/B5dRXT8Xerw}{N=3} &
{22}
\end{tabular}
\caption{Quick-reference to videos mentioned in the paper. Column ``PL\#'' is the entry within the Youtube playlist \cite{dsr_math_intell_playlist}.}
\label{tab:playlist}
\end{table}

\noindent \textbf{Postscript}: Elegant proofs for some of the $N>3$ invariants have recently appeared \cite{akopyan2020-invariants,bialy2020-invariants}.

\begin{acknowledgements}
We would like to thank Sergei Tabachnikov, Richard Schwartz, Arseniy Akopyan, Olga Romaskevich, and Igor Minevich for invaluable mentoring, and key proofs. We would like to thank Jorge Zubelli, Marcos Craizer, Ethan Cotterill, and Matt Perlmutter for inviting us for talks and Paulo Ney de Souza for his encouragement, mathematical, and editorial help. We also thank Mark Helman for proofs for the ellipticity of new triangular centers and Dominique Laurain for a surprising expression for the radius of the Cosine Circle.

The second author is fellow of CNPq and coordinator of Project PRONEX/CNPq/FAPEG 2017 10    26 7000 508.
\end{acknowledgements}

\bibliographystyle{spmpsci} 
\bibliography{elliptic_billiards_v2a}

\begin{thebibliography}{10}
\providecommand{\url}[1]{{#1}}
\providecommand{\urlprefix}{URL }
\expandafter\ifx\csname urlstyle\endcsname\relax
  \providecommand{\doi}[1]{DOI~\discretionary{}{}{}#1}\else
  \providecommand{\doi}{DOI~\discretionary{}{}{}\begingroup
  \urlstyle{rm}\Url}\fi

\bibitem{akopyan19_private_meromorphic}
Akopyan, A.: Proofs of constant cosine product, and area ratio for ${N}>3$.
\newblock Private Communication (August, 2019)

\bibitem{akopyan2020-invariants}
Akopyan, A., Schwartz, R., Tabachnikov, S.: Billiards in ellipses revisited
  (2020).
\newblock \urlprefix\url{https://arxiv.org/abs/2001.02934}

\bibitem{bialy2020-invariants}
Bialy, M., Tabachnikov, S.: {Dan Reznik's} identities and more (2020).
\newblock \urlprefix\url{https://arxiv.org/abs/2001.08469}

\bibitem{birkhoff1927}
Birkhoff, G.: On the periodic motions of dynamical systems.
\newblock Acta Mathematica \textbf{50}(1), 359--379 (1927).
\newblock \doi{10.1007/BF02421325}.
\newblock \urlprefix\url{https://doi.org/10.1007/BF02421325}

\bibitem{connes07}
Connes, A., Zagier, D.: A property of parallelograms inscribed in ellipses.
\newblock The American Mathematical Monthly \textbf{114}(10), 909--914 (2007).
\newblock
  \urlprefix\url{https://people.mpim-bonn.mpg.de/zagier/files/amm/114/fulltext.pdf}

\bibitem{dragovic11}
Dragovi\'{c}, V., Radnovi\'{c}, M.: Poncelet Porisms and Beyond: Integrable
  Billiards, Hyperelliptic Jacobians and Pencils of Quadrics.
\newblock Frontiers in Mathematics. Springer, Basel (2011).
\newblock \urlprefix\url{https://books.google.com.br/books?id=QcOmDAEACAAJ}

\bibitem{corentin19}
Fierobe, C.: On the circumcenters of triangular orbits in elliptic billiard
  (2018).
\newblock \urlprefix\url{https://arxiv.org/pdf/1807.11903.pdf}.
\newblock Arxiv

\bibitem{ronaldo16}
Garcia, R.: Centers of inscribed circles in triangular orbits of an elliptic
  billiard (2016).
\newblock \urlprefix\url{https://arxiv.org/pdf/1607.00179v1.pdf}.
\newblock Arxiv

\bibitem{ronaldo19}
Garcia, R.: Elliptic billiards and ellipses associated to the 3-periodic
  orbits.
\newblock American Mathematical Monthly \textbf{126}(06), 491--504 (2019).
\newblock \urlprefix\url{https://doi.org/10.1080/00029890.2019.1593087}

\bibitem{ronaldo19a}
Garcia, R., Reznik, D., Koiller, J.: New properties of triangular orbits in
  elliptic billiards (2019).
\newblock In preparation

\bibitem{helman19}
Helman, M.: The loci of ${X}_i,i=7,40,57,63,142,144$ are elliptic.
\newblock Private Communication (August, 2019)

\bibitem{himmelstrand12}
Himmelstrand, M., Wilén, V., Saprykina, M.: A survey of dynamical billiards
  (2012).
\newblock \urlprefix\url{https://bit.ly/2kfYHkC}

\bibitem{johnson29}
Johnson, R.A.: Modern Geometry: An Elementary Treatise on the Geometry of the
  Triangle and the Circle.
\newblock Houghton Mifflin, Boston, MA (1929)

\bibitem{jovanovic11}
Jovanović, B.: What are completely integrable hamilton systems.
\newblock The Teaching of Mathematics \textbf{13}(1), 1--14 (2011)

\bibitem{etc}
Kimberling, C.: Encyclopedia of triangle centers (2019).
\newblock
  \urlprefix\url{https://faculty.evansville.edu/ck6/encyclopedia/ETC.html}

\bibitem{dominique19}
Laurain, D.: Formula for the radius of the orbits' excentral cosine circle.
\newblock Private Communication (August, 2019)

\bibitem{minevich17}
Minevich, I., Morton, P.: Synthetic foundations of cevian geometry, {III}: The
  generalized orthocenter.
\newblock Journal of Geometry \textbf{108}, 437--455 (2017).
\newblock \doi{10.1007/s00022-016-0350-2}.
\newblock \urlprefix\url{https://arxiv.org/pdf/1506.06253.pdf}

\bibitem{dsr_applet_x12345}
Reznik, D.: Applet showing the locus of several triangular centers (2019).
\newblock \urlprefix\url{https://editor.p5js.org/dreznik/full/i1Lin7lt7}

\bibitem{dsr_math_intell_playlist}
Reznik, D.: {YouTube} playlist for mathematical intelligencer (2019).
\newblock \urlprefix\url{https://bit.ly/2kTvPPr}

\bibitem{reznik_media}
Reznik, D., Garcia, R., Koiller, J.: Media for elliptic billards and family of
  orbits (2019).
\newblock
  \urlprefix\url{https://dan-reznik.github.io/Elliptical-Billiards-Triangular-Orbits/videos.html}

\bibitem{reznik_web}
Reznik, D., Garcia, R., Koiller, J.: New properties of triangular orbits in
  elliptic billiards (website) (2019).
\newblock
  \urlprefix\url{https://dan-reznik.github.io/Elliptical-Billiards-Triangular-Orbits/}

\bibitem{reznik2020-loci}
Reznik, D., Garcia, R., Koiller, J.: Loci of triangular orbits in elliptic
  billiards  (2020).
\newblock To appear

\bibitem{olga14}
Romaskevich, O.: On the incenters of triangular orbits on elliptic billiards.
\newblock Enseign. Math. \textbf{60}(3-4), 247--255 (2014).
\newblock \doi{10.4171/LEM/60-3/4-2}.
\newblock \urlprefix\url{https://arxiv.org/pdf/1304.7588.pdf}

\bibitem{olga19_mitten}
Romaskevich, O.: Proof the mittenpunkt is stationary.
\newblock Private Communication (May, 2019)

\bibitem{sergei2016}
Schwartz, R., Tabachnikov, S.: Centers of mass of {P}oncelet polygons, 200
  years after.
\newblock Math. Intelligencer \textbf{38}(2), 29--34 (2016).
\newblock \doi{10.1007/s00283-016-9622-9}.
\newblock \urlprefix\url{http://www.math.psu.edu/tabachni/prints/Poncelet5.pdf}

\bibitem{sergei2002}
Tabachnikov, S.: Ellipsoids, complete integrability and hyperbolic geometry.
\newblock Moscow Mathematical Journal \textbf{2}(1), 185--198 (2002).
\newblock \urlprefix\url{http://ftp.math.psu.edu/tabachni/prints/integr.pdf}

\bibitem{sergei91}
Tabachnikov, S.: Geometry and Billiards, \emph{Student Mathematical Library},
  vol.~30.
\newblock American Mathematical Society, Providence, RI (2005).
\newblock \doi{10.1090/stml/030}.
\newblock
  \urlprefix\url{http://www.personal.psu.edu/sot2/books/billiardsgeometry.pdf}.
\newblock Mathematics Advanced Study Semesters, University Park, PA

\bibitem{sergei2016proj}
Tabachnikov, S.: Projective configuration theorems: old wine into new
  wineskins.
\newblock In: S.~Dani, A.~Papadopoulos (eds.) Geometry in History, pp.
  401--434. Springer Verlag (2019).
\newblock \urlprefix\url{https://arxiv.org/pdf/1607.04758.pdf}

\bibitem{sergei19_private_cosine_sum_expression}
Tabachnikov, S.: Proof of expression for constant cosine sum for any $n$.
\newblock Private Communication (October, 2019)

\bibitem{sergei19_private_circles}
Tabachnikov, S.: Proofs of stationary circle radius for ${N}>3$.
\newblock Private Communication (July, 2019)

\bibitem{sergei19_private_meromorphic}
Tabachnikov, S., Schwartz, R.: Proof of constant cosine sum ${N}>3$.
\newblock Private Communication (July, 2019)

\bibitem{mw}
Weisstein, E.: Mathworld (2019).
\newblock \urlprefix\url{http://mathworld.wolfram.com}

\end{thebibliography}

\end{document}